\documentclass[11pt]{amsart}
\usepackage[latin1]{inputenc}
\usepackage{amsfonts}
\usepackage{amsthm}
\usepackage{amsmath}
\usepackage{amssymb}
\usepackage{amscd}
\usepackage{verbatim}
\usepackage{multicol}
\usepackage{tabularx}

\theoremstyle{definition}
\newtheorem{thm}{Theorem}[section]
\newtheorem{prop}[thm]{Proposition}
\newtheorem{cor}[thm]{Corollary}
\newtheorem{con}[thm]{Conjecture}
\newtheorem{lem}[thm]{Lemma}
\newtheorem{defn}[thm]{Definition}

\newtheorem*{rem}{Remark}

\numberwithin{equation}{section}

\def\blb#1{\text{$\mathbb{#1}$}}
\def\cal#1{\text{$\mathcal{#1}$}}

\def\ord#1^#2{#1$^{\text{#2}}$}

\def\lie#1{\mathfrak{#1}}
\def\tlie#1{\tilde{\mathfrak{#1}}}
\def\hlie#1{\hat{\mathfrak{#1}}}

\def\uqr#1^#2{\text{$U_q^{#2}(\lie #1)$}}

\def\uqhr#1^#2{\text{$U_q^{#2}(\hlie #1)$}}
\def\us#1^#2{\text{$U_{\xi}^{#2}(\lie #1)$}}
\def\ush#1^#2{\text{$U_{\xi}^{#2}(\hlie #1)$}}
\def\dus#1^#2{\text{$\dot{U}_{\xi}^{#2}(\lie #1)$}}
\def\dush#1^#2{\text{$\dot{U}_{\xi}^{#2}(\hlie #1)$}}
\def\gb#1{{\mbox{\boldmath $#1$}}}
\def\gbr#1{{\mbox{\boldmath ${\rm #1}$}}}
\def\wtl{{\rm wt}_\ell}
\def\wt{{\rm wt}}

\def\gr{{\rm gr}}
\def\chl{{\rm char}_\ell}
\def\ch{{\rm char}}

\def\opl_#1^#2{\text{\scriptsize$\bigoplus\limits_{\text{\footnotesize$#1$}}^{\text{\footnotesize$#2$}}$}}
\def\otm_#1^#2{\text{\scriptsize$\bigotimes\limits_{\text{\footnotesize$#1$}}^{\text{\footnotesize$#2$}}$}}
\def\wcal#1{{\mbox{$\widetilde{\cal #1}$}}}

\def\bgb#1{{\mbox{$\overline{\gb #1}$}}}

\renewcommand{\thefootnote}

\begin{document}

\title[]{Restricted limits of minimal affinizations}
\author{Adriano Moura}
\thanks{}
\address{UNICAMP - IMECC, Campinas - SP - Brazil, 13083-970.}
\email{aamoura@ime.unicamp.br}

\maketitle
\centerline{\small{
\begin{minipage}{350pt}
{\bf Abstract:} We obtain character formulas of minimal affinizations of representations of quantum groups when the underlying simple Lie algebra is orthogonal and the support of the highest weight is contained in the first three nodes of the Dynkin diagram. We also give a framework for extending our techniques to a more general situation. In particular, for the orthogonal algebras and a highest weight supported in at most one spin node, we realize the restricted classical limit of the corresponding minimal affinizations as a quotient of a module given by generators and relations and, furthermore, show that it projects onto the submodule generated by the top weight space of the tensor product of appropriate restricted Kirillov-Reshetikhin modules. We also prove a conjecture of Chari and Pressley regarding the equivalence of certain minimal affinizations in type $D_4$.
\end{minipage}}}

\setcounter{section}{0}

\section*{Introduction}

The representation theory of  affine Kac-Moody algebras and their quantum groups has been intensively studied from a broad range of perspectives in the last two decades. In this paper we focus on non-twisted quantum affine Kac-Moody algebras and their finite-dimensional representations. Let $\lie g$ be a finite-dimensional simple Lie algebra over the complex numbers, $\tlie g=\lie g\otimes \mathbb C[t,t^{-1}]$ the associated loop algebra, and $U_q(\lie g), U_q(\tlie g)$ their Drinfeld-Jimbo quantum groups over $\mathbb C(q)$, where $q$ is an indeterminate. The affine Kac-Moody algebra is a one-dimensional central extension of $\tlie g$ but, since the center acts trivially on finite-dimensional modules, it suffices to consider the loop algebra. It turns out that the finite-dimensional representations of $U_q(\tlie g)$ are $\ell$-weight modules, i.e., every vector is a linear combination of common generalized eigenvectors for $U_q(\tlie h)$ where $\lie h$ is a fixed Cartan subalgebra of $\lie g$ and $\tlie h=\lie h\otimes \mathbb C[t,t^{-1}]$.  Moreover, the simple modules are highest-$\ell$-weight and the set of all dominant $\ell$-weights is in bijection with the monoid $\cal P_q^+$ of $n$-tuples of polynomials in one variable with constant term $1$, where $n$ is the rank of  $\lie g$. The set of all $\ell$-weights corresponds to the group $\cal P_q$ associated to $\cal P_q^+$. By specializing $q$ at $1$ one recovers the finite-dimensional representation theory of $\tlie g$.

Given a nonzero complex number $a$, let ${\rm ev}_a:\tlie g\to \lie g$ be the evaluation map $x\otimes f(t)\mapsto f(a)x$. If $V$ is a $\lie g$-module, one can consider the pullback $V(a)$ of $V$ by ${\rm ev}_a$. In particular, every irreducible $\lie g$-module can be turned into a $\tlie g$-module.
In the quantum case, unless $\lie g$ is of type $A$, there is no analogue of the evaluation map and, in fact, most often, an irreducible $U_q(\lie g)$-module cannot be turned into a $U_q(\tlie g)$-module. By allowing the underlying vector space to be enlarged in a ``controlled'' way, a concept of quantum affinization of an irreducible $U_q(\lie g)$-module was introduced in \cite{cha:minr2}.  Two affinizations are said to be equivalent if they have isomorphic $U_q(\lie g)$-structures. It follows from the classification of the finite-dimensional irreducible $U_q(\tlie g)$-modules that every finite-dimensional irreducible $U_q(\lie g)$-module has at least one equivalence class of affinizations. Moreover, there are finitely many equivalence classes of affinizations  and the usual partial order on the weight lattice $P$ of $\lie g$ induces a partial order on the set of equivalence classes of affinizations of a given irreducible $U_q(\lie g)$-module.  Representatives of the minimal elements with respect to this partial order are called minimal affinizations. Although an almost complete classification of the highest $\ell$-weights of equivalence classes of minimal affinizations was obtained in \cite{cha:minr2,cp:minsl,cp:minnsl,cp:minirr}, their structure remained essentially unknown except when $\lie g$ is of type $A$ or $B_2$. Further progress was made after the introduction of the concept of $q$-characters in \cite{freer:qchar}, which we prefer to call $\ell$-characters as explained in Section \ref{ss:qchar}.

The $\ell$-character of a finite-dimensional $U_q(\tlie g)$-module $V$ is the associated element $\chl(V)$ of the integral group ring $\mathbb Z[\cal P_q]$ which records the dimensions of the $\ell$-weight spaces of $V$. Given $\gb\lambda\in\cal P_q^+$, let us denote by $V_q(\gb\lambda)$ the irreducible $U_q(\tlie g)$-module with highest $\ell$-weight $\gb\lambda$. Finding formulas for the $\ell$-character of $V_q(\gb\lambda)$ is still an an open problem in general. In \cite{freem:qchar}, E. Frenkel and E. Mukhin defined an algorithm, now widely known as the Frenkel-Mukhin algorithm which for agiven $\gb\lambda\in \cal P_q^+$ returns an element of $\mathbb Z[\cal P_q]$ that was  conjectured to be $\chl(V_q(\gb\lambda))$. The conjecture was proved for certain situations in \cite{freem:qchar}, but it has recently been shown that this is not always the case \cite{nanak:fma}. However, even in the situations for which the conjecture holds, the task of translating the information given by the algorithm  into general closed formulas remains a challenge.  For further details on the theory of $\ell$-characters, beside the aforementioned literature,  we refer the reader to the very recent survey \cite{chhe:beyond} and the references therein. We remark that in \cite{nanak:trd,nanak:trc} the authors give path-tableaux descriptions of Jacobi{-}Trudi determinants which, conjecturally, coincide with the $\ell$-characters if $\lie g$ is of classical type. This conjecture has been partially proved if $\lie g$ is of type $B$ in \cite{her:min}(see also \cite{chhe:beyond}).

Another approach for studying minimal affinizations is by considering their classical limit. Even though most of the $\ell$-character information is lost in this process, it provides an effective tool to study their $U_q(\lie g)$-structure, i.e., their characters. The $U_q(\lie g)$-structure of the minimal affinizations belonging to the family of Kirillov-Reshetikhin modules was obtained in \cite{cha:fer} partially using this approach. The proof consisted in showing that the conjectural character was both a lower and an upper bound for the character of the given Kirillov-Reshetikhin module. While the latter was proved by working with the classical limit, the proof of the former was done in the quantum context.  Later on, it was shown in \cite{cm:kr,cm:krg} that both ``upper and lower bound'' parts of the proofs of the results of \cite{cha:fer} could be performed by working with the current algebra $\lie g[t]=\lie g\otimes\mathbb C[t]$. Characters of Kirillov-Reshetikhin modules for twisted affine algebras were also obtained in \cite{cm:kr,cm:krg} in this manner. The Kirillov-Reshetikhin modules were introduced in \cite{kr:krm} (in the context of Yangians rather than quantumm affine algebras) in connection with the Bethe Ansatz. They are the minimal affinizations of the irreducible $U_q(\lie g)$-modules whose highest weights are multiples of  the fundamental weights of $\lie g$.

The main goal of the present paper is to initiate a program for extending the approach of \cite{cm:kr,cm:krg} to more general minimal affinizations other than Kirillov-Reshetikhin modules. We prove several partial results in this direction and carry out the whole program in the simplest cases. In particular, we obtain character formulas for minimal affinizations in the case that $\lie g$ is orthogonal and the support of the highest weigh is contained in the first three nodes of the Dynkin diagram of $\lie g$. We now give a summary of our results.

Given a dominant integral weight $\lambda=\sum m_i\omega_i$ (where $\omega_i, i=1,\dots,n$, are the fundamental weights of $\lie g$), we define restricted graded $\lie g[t]$-modules $M(\lambda)$ and $T(\lambda)$. The former is given by generators and relations while the latter is the submodule generated by the top weight space of $\otimes_i M(m_i\omega_i)$. We conjecture that these modules are isomorphic.  This is a generalization of one of the main results of \cite{cm:kr,cm:krg}. The conjecture clearly holds for type $A$. The defining relations for the module $M(\lambda)$ are, roughly speaking, the intersection of the relations satisfied by the corresponding restricted Kirillov-Reshetikhin modules $M(m_i\omega_i)$. In particular, it is immediate that $T(\lambda)$ is a quotient of $M(\lambda)$. We prove this conjecture when $\lie g$ is orthogonal and $\lambda$ is supported only in the first three nodes of the Dynkin diagram of $\lie g$. If $\lie g$ is of type $D$, the proof also works in the case that both spin nodes are in the support of $\lambda$. As a byproduct of the proof, we obtain the characters of the modules $M(\lambda)$ in these cases.
Namely, assume $\lie g$ is of type $B_n$ and that the nodes of the Dynkin diagram of $\lie g$ are labeled as in \cite{hum:book}. Given $\lambda=m_1\omega_1+m_2\omega_2+m_3\omega_3$, consider the set $\cal A=\{\gbr r=(r_1,r_2,r_3)\in\mathbb Z_{\ge 0}^3: r_1+r_2\le [a_3m_3], r_2\le m_1, r_3\le [a_2m_2]\}$, where $[m]$ denotes the integer part of the rational number $m$, $a_n=1/2$, and $a_i=1$ for $i\ne n$. Then, we have an isomorphism of $\lie g$-modules:
\begin{gather*}\tag{1}
M(\lambda)\cong \opl_{\gbr r\in\cal A}^{} V((m_1+r_1-r_2)\omega_1+(m_2+r_2-a_2^{-1}r_3)\omega_2+(m_3-a_3^{-1}(r_1+r_2))\omega_3).
\end{gather*}
Here, $V(\mu)$ denotes the irreducible $\lie g$-module of highest weight $\mu\in P^+$. If $\lie g$ is of type $D_n$ with $n\ge 5$ and $\lambda=m_1\omega_1+m_2\omega_2+m_3\omega_3+m_{n-1}\omega_{n-1}+m_n\omega_n$, the $\lie g$-structure of $M(\lambda)$ is given by (1) as well (in this case $a_i=1$ for all $i$). If $n=4$ and $\lambda\in P^+$, then $M(\lambda)\cong\opl_{r=0}^{m_2} V(\lambda-r\omega_2)$ as a $\lie g$-module.

On the other hand, by regarding the classical limit of a minimal affinization $V_q(\gb\lambda)$ as a $\lie g[t]$-module and then shifting the associated spectral parameter to zero, we obtain modules $L(\gb\lambda)$ which we call the restricted limit of $V_q(\gb\lambda)$. Let $\lambda$ be the maximal weight of $V_q(\gb\lambda)$. We prove that $T(\lambda)$ is a quotient of $L(\gb\lambda)$ (Proposition \ref{p:TqM}). Moreover, for orthogonal $\lie g$, we prove that $L(\gb\lambda)$ is a quotient of $M(\lambda)$ provided that $\lambda$ is supported in a connected subdiagram of type $A$ if $\lie g$ is of type $D$ (Proposition \ref{p:MqN}). Therefore, if indeed $M(\lambda)$ is isomorphic to $T(\lambda)$ as conjectured, it would follow that they are also isomorphic to $L(\gb\lambda)$ in the above cases.  In particular, equation (1) above describes the $U_q(\lie g)$-structure of $V_q(\gb\lambda)$ when $\lie g$ is orthogonal and $\lambda$ is supported only on the first three nodes of the Dynkin diagram of $\lie g$ (and possibly on one of the spin nodes if $\lie g$ is of type $D$). For $\lie g$ of type $B_2$, the same result was obtained in \cite{cha:minr2} by working purely in the quantum setting. If $\lie g$ is of type $B_n$ and the value of $\lambda$ on the coroot associated to the spin node is even, then the $\ell$-character (and hence the character) of $V_q(\gb\lambda)$ can be computed using the tableaux expression of Jacobi-Trudi determinant (see \cite[\S 7.6]{chhe:beyond}).  We expect that, if the minimal connected subdiagram of the Dynkin diagram of $\lie g$ containing the support of $\lambda$ does not contain a subdiagram of type $D_4$ (in which case $V_q(\lambda)$ has a unique equivalence class of minimal affinizations),  Proposition \ref{p:MqN} remains valid and, hence, that the modules $T(\lambda), M(\lambda)$, and $L(\gb\lambda)$ are isomorphic. To keep the length of the present paper within reasonable limits, we leave the quest of pursuing the proofs of these conjectures in a more general setting to a forthcoming publication.

When $V_q(\lambda)$ has more than one equivalence class of minimal affinizatios, it is certainly not true that $L(\gb\lambda)$ is a quotient of $M(\lambda)$ (in fact, it is the other way round). It was proved in \cite{cp:minsl} that, if $\lambda$ is supported in the triply connected node of the Dynkin diagram of $\lie g$, then there are exactly three equivalence classes of minimal affinizations. We define $\lie g[t]$-modules $M_k(\lambda), k=1,2,3$, and prove that $L(\gb\lambda)$ is a quotient of $M_k(\lambda)$ for exactly one value of $k$. Naturally, we expect that $L(\gb\lambda)$ is isomorphic to the appropriate $M_k(\lambda)$. We prove that this is so if $\lie g$ is of type $D_4$ and obtain the character of $M_k(\lambda)$ in this case. Namely, let $\lambda=m_1\omega_1+m_2\omega_2+m_3\omega_3+m_4\omega_4$, where the triply connected node is labeled by $4$, suppose $\{i,j,k\}=\{1,2,3\}$, and set $\cal A_k = \{\gbr r\in\mathbb Z_{\ge 0}^3: r_1\le m_k, r_1+r_2\le \min\{m_i,m_j\}, r_3\le m_4\}$. Then, we have an isomorphism of $\lie g$-modules
\begin{gather*}\tag{2}
M_k(\lambda)\cong\opl_{\gbr r\in\cal A_k}^{} V(\lambda-(r_1-r_2)\omega_k - (r_1+r_2)(\omega_i+\omega_j)-(r_3-r_1)\omega_4).
\end{gather*}

If $\lambda$ is not supported in the triply connected node, it was proved in \cite[Theorem 2.2]{cp:minirr} that the number of equivalence classes of minimal affinizations of $V_q(\lambda)$ grows as $\lambda$ ``grows''. Although we do not have a general conjecture in this case yet,  the definition of $M_k(\lambda)$ makes sense in this case as well and its character is computed in the same way as in the previous case.  Moreover, the same proof we applied to the previous case in type $D_4$ also proves that, if $\gb\lambda$ satisfies the conditions (a)$_{i,j}$ or (b)$_{i,j}$ of \cite[Theorem 2.2]{cp:minirr}, then $L(\gb\lambda)$ is isomorphic to  $M_k(\lambda)$ for the appropriate value of $k$ and its character is given by equation (2) above. In particular, this proves the conjecture of \cite{cp:minirr} saying that the modules $V_q(\gb\lambda)$ with $\gb\lambda$ satisfying conditions (a)$_{i,j}$ of \cite[Theorem 2.2]{cp:minirr} are equivalent to those with $\gb\lambda$ satisfying conditions (b)$_{i,j}$ of that same theorem.

The techniques employed to prove Propositions \ref{p:TqM} and \ref{p:MqN} (and their analogues in the case of multiple equivalence classes of minimal affinizations)  make use of the results of \cite{cha:braid} in an essential way. Moreover, for the proof of Proposition \ref{p:MqN} we also use  some partial information on $\ell$-characters by combining the Frenkel-Mukhin algorithm  with results proved in \cite{cm:chb,freem:qchar,her:min}. The paper is organized as follows. In Sections \ref{s:algs} and \ref{s:fdimreps} we review some structural results of the algebras $\lie g,\tlie g$, and their quantum counterparts as well as some basic results of the finite-dimensional representation theory of these algebras. In Section \ref{s:min}, after reviewing the partial classification of minimal affinizations, we define the modules $M(\lambda), T(\lambda)$, and $L(\gb\lambda)$, and state our main results and conjectures regarding them. The proofs are given in Sections \ref{s:tensor} and \ref{s:char}. The case of multiple equivalence classes of minimal affinizations is treated in Subsections \ref{ss:multmar} and \ref{ss:multmai}.

\vskip15pt
\noindent{\bf Acknowledgements:}  This work was partially supported by CNPq and FAPESP. The author thanks D. Jakeli\'c for helpful discussions.

\section{Quantum and classical loop algebras}\label{s:algs}

Throughout the paper, let $\mathbb C, \mathbb R,\mathbb Z,\mathbb Z_{\ge m}$ denote the sets of complex numbers, reals, integers,  and integers bigger or equal $m$, respectively. Given a ring $\mathbb A$, the underlying multiplicative group of units is denoted by $\mathbb A^\times$. The dual of a vector space $V$ is denoted by $V^*$. The symbol $\cong$ means ``isomorphic to''. 

\subsection{Classical algebras}\label{ss:clalg}

Let $I=\{1,\dots,n\}$ be the set of vertices of a finite-type connected Dynkin diagram labeled as in \cite{hum:book}
and let $\lie g$ be the associated simple Lie algebra over $\mathbb C$ with a fixed Cartan subalgebra $\lie h$. Fix a set of positive roots $R^+$ and let
$$\lie n^\pm = \opl_{\alpha\in R^+}^{} \lie g_{\pm\alpha} \quad\text{where}\quad \lie g_{\pm\alpha} = \{x\in\lie g: [h,x]=\pm\alpha(h)x, \ \forall \ h\in\lie h\}.$$
The simple roots will be denoted by $\alpha_i$, the fundamental weights by $\omega_i$, while $Q,P,Q^+,P^+$ will denote the root and weight lattices with corresponding positive cones, respectively.
Let also $h_i\in\lie h$, be the co-root associated to $\alpha_i, i\in I$. We equip $\lie h^*$ with the partial order $\lambda\le \mu$ iff $\mu-\lambda\in Q^+$.
We denote by $\cal W$ the Weyl group of $\lie g$
and let $w_0$ be the longest element of $\cal W$. Given $\lambda\in P$ set
\begin{equation}
\lambda^*=-w_0\lambda.
\end{equation}
Recall that, if $\lambda\in P^+$, then $\lambda^*\in P^+$ as well. Let $C = (c_{ij})_{i,j\in I}$ be the Cartan matrix of $\lie g$, i.e., $c_{ij}=\alpha_j(h_i)$, and let $D = {\rm diag}(d_i:i\in I)$ be such that the numbers $d_i$ are coprime positive integers satisfying $DC$ is symmetric.

The subalgebras $\lie g_{\pm\alpha}, \alpha\in R^+$, are one-dimensional and $[\lie g_{\pm\alpha},\lie g_{\pm\beta}]=\lie g_{\pm\alpha\pm\beta}$ for every $\alpha,\beta\in R^+$. We denote by $x^\pm_\alpha$ any generator of $\lie g_{\pm\alpha}$. In particular, if $\alpha+\beta\in R^+$, $[x^\pm_\alpha,x^\pm_\beta]$ is a nonzero generator of $\lie g_{\pm\alpha\pm\beta}$ and we simply write  $[x^\pm_\alpha,x^\pm_\beta]=x_{\alpha+\beta}^\pm$. For each subset $J$ of $I$ let $\lie g_J$ be the Lie subalgebra of $\lie g$ generated by $x_{\alpha_j}^\pm, j\in J$, and define $\lie n^\pm_J, \lie h_J$ in the obvious way. Let also $Q_J$ be the subgroup of $Q$ generated by $\alpha_j, j\in J$, and $R^+_J=R^+\cap Q_J$. Given $\lambda\in P$, let $\lambda_J$ be the restriction of $\lambda$ to $\lie h_J^*$ and $\lambda^J\in P$ be such that $\lambda^J(h_j)=\lambda(h_j)$ if $j\in J$ and $\lambda^J(h_j)=0$ otherwise. By abuse of language, we will refer to any subset $J$ of $I$ as a subdiagram of the Dynkin diagram of $\lie g$. The support of $\mu\in P$ is defined to be the subdiagram ${\rm supp}(\mu)\subseteq I$ given by ${\rm supp}(\mu)=\{i\in I:\mu(h_i)\ne 0\}$. Let also $\overline{\rm supp}(\mu)$ be the minimal connected subdiagram of $I$ containing ${\rm supp}(\mu)$.

If  $\lie a$ is a Lie algebra over $\mathbb C$, define its loop algebra to be $\tlie a=\lie a\otimes_{\mathbb C}  \mathbb C[t,t^{-1}]$ with bracket given by $[x \otimes t^r,y \otimes t^s]=[x,y] \otimes t^{r+s}$. Clearly $\lie a\otimes 1$ is a subalgebra of $\tlie a$ isomorphic to $\lie a$ and, by abuse of notation, we will continue denoting its elements by $x$ instead of $x\otimes 1$. We also consider the current algebra $\lie a[t]$ which is the subalgebra of $\tlie a$ given by $\lie a[t]=\lie a\otimes \mathbb C[t]$. Then $\tlie g = \tlie n^-\oplus \tlie h\oplus \tlie n^+$ and $\tlie h$ is an abelian subalgebra and similarly for $\lie g[t]$. The elements $x_\alpha^\pm\otimes t^r$ and $h_i\otimes t^r$ will be denoted by $x_{\alpha,r}^\pm$ and $h_{\alpha_i,r}$, respectively. Diagram subalgebras $\tlie g_J$ are defined in the obvious way.

Let $U(\lie a)$ denote the universal enveloping algebra of a Lie algebra $\lie a$. Then $U(\lie a)$ is a subalgebra of $U(\tlie a)$ and multiplication establishes isomorphisms of vector spaces
$$U(\lie g)\cong U(\lie n^-)\otimes U(\lie h)\otimes U(\lie n^+) \qquad\text{and}\qquad U(\tlie g)\cong U(\tlie n^-)\otimes U(\tlie h)\otimes U(\tlie n^+).$$
The assignments $\triangle: \lie a\to U(\lie a)\otimes U(\lie a), x\mapsto x\otimes 1+1\otimes x$, $S:\lie a\to \lie a,
x\mapsto -x$,  and $\epsilon: \lie a\to \mathbb C, x\mapsto 0$, can
be uniquely extended so that $U(\lie a)$ becomes a Hopf algebra with
comultiplication $\triangle$, antipode $S$, and  counit $\epsilon$.

Given $a\in\mathbb C$, let $\tau_a$ be the Lie algebra automorphism of $\lie a[t]$ defined by $\tau_a(x\otimes f(t))=x\otimes f(t-a)$ for every $x\in\lie a$ and every $f(t)\in\mathbb C[t]$. If $a\ne 0$, let ${\rm ev}_a:\tlie a\to\lie a$ be the evaluation map $x\otimes f(t)\mapsto f(a)x$. We also denote by $\tau_a$ and ${\rm ev}_a$ the induced maps $U(\lie a[t])\to U(\lie a[t])$ and $U(\tlie a)\to U(\lie a)$, respectively.

For each $i\in I$ and  $r\in\mathbb Z$, define elements $\Lambda_{i,r}\in U(\tlie h)$ by the following equality of formal power series in the variable $u$:
\begin{equation}\label{e:Lambdadef}
\sum_{r=0}^\infty \Lambda_{i,\pm r} u^r = \exp\left( - \sum_{s=1}^\infty \frac{h_{\alpha_i,\pm s}}{s} u^s\right).
\end{equation}

\subsection{Quantum algebras}

Let $\mathbb C(q)$ be the ring of rational functions on an indeterminate $q$ and $\mathbb A=\mathbb C[q,q^{-1}]$. Given $p=q^{k}$ for some $k\in\mathbb Z\backslash\{0\}$, define
\begin{equation*}
[m]_p=\frac{p^m -p^{-m}}{p -p^{-1}},\ \ \ \ [m]_p!
=[m]_p[m-1]_p\ldots [2]_p[1]_p,\ \ \ \ \left[\begin{matrix} m\\
r\end{matrix}\right]_p = \frac{[m]_p!}{[r]_p![m-r]_p!},
\end{equation*}
for $r,m\in\mathbb Z_{\ge 0}$, $m\ge r$. Notice that $\left[\begin{matrix} m\\
r\end{matrix}\right]_p\in\mathbb A$.

Set $q_i=q^{d_i}$. The quantum loop algebra  $\
U_q(\tlie g)$ of $\lie g$ is  the algebra with
generators $x_{i,r}^{{}\pm{}}$ ($i\in I$, $r\in\blb Z$), $k_i^{{}\pm
1}$ ($i\in I$), $h_{i,r}$ ($i\in I$, $r\in \blb Z\backslash\{0\}$)
and the following defining relations:
\begin{align*}
k_ik_i^{-1} = k_i^{-1}k_i& =1, \ \  k_ik_j =k_jk_i,\\
k_ih_{j,r}& =h_{j,r}k_i,\\
k_ix_{j,r}^\pm k_i^{-1} &= q_i^{{}\pm c_{ij}}x_{j,r}^{{}\pm{}},\ \ \\
[h_{i,r},h_{j,s}]=0,\; \; & [h_{i,r} , x_{j,s}^{{}\pm{}}] = \pm\frac1r[rc_{ij}]_{q^i}x_{j,r+s}^{{}\pm{}},\\
x_{i,r+1}^{{}\pm{}}x_{j,s}^{{}\pm{}} -q_i^{{}\pm c_{ij}}x_{j,s}^{{}\pm{}}x_{i,r+1}^{{}\pm{}} &
=q_i^{{}\pm     c_{ij}}x_{i,r}^{{}\pm{}}x_{j,s+1}^{{}\pm{}} -x_{j,s+1}^{{}\pm{}}x_{i,r}^{{}\pm{}},\\
[x_{i,r}^+ , x_{j,s}^-]=\delta_{i,j} & \frac{ \psi_{i,r+s}^+ - \psi_{i,r+s}^-}{q_i - q_i^{-1}},\\
\sum_{\sigma\in S_m}\sum_{k=0}^m(-1)^k
\left[\begin{matrix}m\\k\end{matrix}\right]_{q_i}
x_{i, r_{\sigma(1)}}^{{}\pm{}}\ldots x_{i,r_{\sigma(k)}}^{{}\pm{}} &
x_{j,s}^{{}\pm{}} x_{i, r_{\sigma(k+1)}}^{{}\pm{}}\ldots x_{i,r_{\sigma(m)}}^{{}\pm{}} =0,\ \ \text{if $i\ne j$},
\end{align*}
for all sequences of integers $r_1,\ldots, r_m$, where $m
=1-c_{ij}$, $S_m$ is the symmetric group on $m$ letters, and
the $\psi_{i,r}^{{}\pm{}}$ are determined by equating powers of
$u$ in the formal power series
$$\Psi_i^\pm(u) = \sum_{r=0}^{\infty}\psi_{i,\pm r}^{\pm}u^r = k_i^{\pm 1} \exp\left(\pm(q_i-q_i^{-1})\sum_{s=1}^{\infty}h_{i,\pm s} u^s\right).$$

Denote by $U_q(\tlie n^\pm), U_q(\tlie h)$  the subalgebras of $U_q(\tlie g)$ generated by $\{x_{i,r}^\pm\}, \{k_i^{\pm1}, h_{i,s}\}$, respectively. Let  $U_q(\lie g)$  be  the subalgebra generated by $x_i^\pm:=x_{i,0}^\pm, k_i^{\pm 1}, i\in I,$ and define $U_q(\lie n^\pm), U_q(\lie h)$ in the obvious way.  $U_q(\lie g)$ is a subalgebra of $U_q(\tlie g)$ and multiplication establishes isomorphisms of $\mathbb C(q)$-vectors spaces:
$$U_q(\lie g) \cong U_q(\lie n^-) \otimes U_q(\lie h) \otimes U_q(\lie n^+) \qquad\text{and}\qquad U_q(\tlie g) \cong U_q(\tlie n^-) \otimes U_q(\tlie h) \otimes U_q(\tlie n^+).$$

Let $J\subseteq I$ and consider the subalgebra $U_q(\tlie g_J)$ generated by $k_j^{\pm 1}, h_{j,r}, x^\pm_{j,s}$ for all $j\in J, r,s\in \blb Z, r\ne 0$. If $J=\{j\}$, the algebra $U_q(\tlie g_j):=U_q(\tlie g_J)$ is isomorphic to $U_{q_j}(\tlie{sl}_2)$. Similarly we define the subalgebra $U_q(\lie g_J)$, etc.

For $i\in I, r\in \mathbb Z, k\in\mathbb Z_{\ge 0}$, define $(x_{i,r}^\pm)^{(k)} = \frac{(x_{i,r}^\pm)^k}{[k]_{q_i}!}$. Define also elements $\Lambda_{i,r}, i\in I, r\in\mathbb Z$ by
\begin{equation}\label{e:Lambdad}
\sum_{r=0}^\infty \Lambda_{i,\pm r} u^{r}=
\exp\left(-\sum_{s=1}^\infty\frac{h_{i,\pm s}}{[s]_{q_i}}u^s\right).
\end{equation}
Note that \begin{equation}
\Psi_i^\pm(u)=k_i^{\pm 1}\frac{\Lambda_i^\pm(q_i^{-1}u)}{\Lambda_i^\pm(q_iu)}
\end{equation}
where the division is that of formal power series in $u$. Although we are denoting the elements $\Lambda_{i,r}$ above by the same symbol as their classical counterparts, this will not create confusion as it will be clear from the context.

Let $U_{\mathbb A}(\tlie g)$ be the $\mathbb A$-subalgebra of $U_q(\tlie g)$ generated by the elements $(x_{i,r}^\pm)^{(k)}, k_i^{\pm 1}$ for $i\in I,r\in\mathbb Z$, and $k\in\mathbb Z_{\ge 0}$. Define $U_{\mathbb A}(\lie g)$ similarly and notice that $U_{\mathbb A}(\lie g)=U_{\mathbb A}(\tlie g)\cap U_q(\lie g)$ . For the proof of the next proposition see \cite[Lemma 2.1]{cha:fer} and the locally cited references.

\begin{prop}
We have $U_q(\tlie g)=\mathbb C(q)\otimes_{\mathbb A} U_{\mathbb A}(\tlie g)$ and $U_q(\lie g)=\mathbb C(q)\otimes_{\mathbb A} U_{\mathbb A}(\lie g)$.\hfill\qedsymbol
\end{prop}

Regard $\mathbb C$ as an $\mathbb A$-module by letting $q$ act as $1$ and set
\begin{equation}
\overline{U_q(\tlie g)} = \mathbb C\otimes_\mathbb A U_\mathbb A(\tlie g) \qquad\text{and}\qquad \overline{U_q(\lie g)} = \mathbb C\otimes_\mathbb A U_\mathbb A(\lie g).
\end{equation}
Denote by $\overline\eta$ the image of $\eta\in U_\mathbb A(\tlie g)$ in $\overline{U_q(\tlie g)}$. The proof of the next proposition can be found in \cite[Proposition 9.2.3]{cp:book} and \cite{lus:book}.

\begin{prop}\label{p:cluq}
$U(\tlie g)$ is isomorphic to the quotient of $\overline{U_q(\tlie g)}$ by the ideal generated by ${\overline k_i}-1$. In particular, the category of $\overline{U_q(\tlie g)}$-modules on which $k_i$ act as the identity operator for all $i\in I$ is equivalent to the category of all $\tlie g$-modules.\hfill\qedsymbol
\end{prop}

The algebra $U_q(\tlie g)$ is a Hopf algebra and induces a Hopf algebra structure (over $\mathbb A$) on $U_\mathbb A(\tlie g)$ (see \cite{cp:book,lus:book}). Moreover, the induced Hopf algebra structure on $U(\tlie g)$ coincides with the usual one.  On $U_q(\lie g)$ we have
\begin{equation}\label{e:comultuqg}
\Delta (x_{i}^+) = x_{i}^+\otimes 1+ k_i\otimes x_{i}^+, \qquad \Delta (x_{i}^-) = x_{i}^-\otimes k_i^{-1}  + 1\otimes x_{i}^-, \qquad \Delta(k_i) = k_i\otimes k_i
\end{equation}
for all $i\in I$.

\begin{lem}\label{l:comcom}
Suppose $x=[x_{i_1}^-,[x_{i_2}^-,\cdots[x_{i_{l-1}}^-,x_{i_l}^-]\cdots]]$. Then $X\in U_\mathbb A(\lie n^-)$ and
$$\Delta(x)\in x\otimes (\prod_{j=1}^l k_{i_j}^{-1}) + 1\otimes x + f(q)y$$
for some $y\in U_\mathbb A(\lie g)\otimes U_\mathbb A(\lie g)$ and some $f(q)\in\mathbb A$ such that $f(1)=0$.
\end{lem}

\begin{proof}
When $l=1$ this is immediate from \eqref{e:comultuqg}. A straightforward induction on $l$ using the relations $k_ix_j^- = q_i^{-c_{i,j}}x_j^-k_i$ completes the proof.
\end{proof}

An expression for the comultiplication $\Delta$ of $U_q(\tlie g)$ in terms of the generators $x^\pm_{i,r}, h_i,r, k_i^{\pm1}$  is not known. The following partial information will suffice for our purposes. Let $X^\pm$ be the subspace of $U_\mathbb A(\tlie n^\pm)$ spanned by $\{x_{j,r}^\pm: j\in I, r\in\mathbb Z\}$.

\begin{lem}\label{l:comult}
$\Delta(x_{i,1}^-)=x_{i,1}^-\otimes k_i + 1\otimes x_{i,1}^- + x$ for some $x\in U_\mathbb A(\lie g)\otimes U_\mathbb A(\lie g)$ such that $\bar x=0$.
\end{lem}

\begin{proof}
It was proved in \cite{be:bgqaa,bcp,dam} (see also Lemma 7.5 of \cite{cp:root}) that $\Delta(x_{i,1}^-)=x_{i,1}^-\otimes k_i + 1\otimes x_{i,1}^- + x$ where $x\in U_\mathbb A(\tlie g)X^-\otimes U_\mathbb A(\tlie g)X^+$. Since the image $\overline{\Delta(x_{i,1}^-)}$ of $\Delta(x_{i,1}^-)$ in $U(\tlie g)\otimes U(\tlie g)$ is $x_{\alpha_i,1}^-\otimes 1 + 1\otimes x_{\alpha_i,1}^-$, it follows that the image of $x$ in $U(\tlie g)\otimes U(\tlie g)$ must be zero.
\end{proof}

The following was also proved in \cite{be:bgqaa,bcp,dam}. Modulo $U_q(\tlie g)X^-\otimes U_q(\tlie g)X^+$, we have
\begin{equation}\label{e:cmLambda}
\Delta(h_{i,r}) =h_{i,s}\otimes 1+1\otimes h_{i,r} \qquad\text{and}\qquad \Delta(\Lambda_{i,r}) = \sum_{s=0}^r \Lambda_{i,r-s}\otimes \Lambda_{i,s} \qquad\text{for all}\qquad r\in\mathbb Z_{\ge 1}.
\end{equation}


We will need the following general result on the dual representation of a tensor product of representations of a Hopf algebra. The proof can be found in \cite{kas:book} for instance.

\begin{prop}\label{p:dualtp}
Let $H$ be a Hopf algebra and $V$ and $W$ be finite-dimensional $H$-modules. Then $(V\otimes W)^*\cong W^*\otimes V^*$.\hfill\qedsymbol
\end{prop}

\subsection{The $\ell$-weight lattice}\label{ss:llattice} Given a field $\mathbb F$ consider the multiplicative group $\cal P_\mathbb F$ of $n$-tuples of rational functions $\gb\mu = (\gb\mu_1(u),\cdots, \gb\mu_n(u))$ with values in $\mathbb F$  such that $\gb\mu_i(0)=1$ for all $i\in I$. We shall often think of $\gb\mu_i(u)$ as a formal power series in $u$ with coefficients in $\mathbb F$. Given $a\in\mathbb F^\times$ and $i\in I$, let $\gb\omega_{i,a}$ be defined by
$$(\gb\omega_{i,a})_j(u) = 1-\delta_{i,j}au.$$
Clearly, if $\mathbb F$ is algebraically closed, $\cal P_\mathbb F$ is the free abelian group generated by these elements which are called fundamental $\ell$-weights. It is also convenient to introduce elements $\gb\omega_{\lambda,a}, \lambda\in P,a\in\mathbb C(q)$, defined by
\begin{equation}
\gb\omega_{\lambda,a} = \prod_{i\in I}(\gb\omega_{i,a})^{\lambda(h_i)}.
\end{equation}
If $\mathbb F$ is algebraically closed, introduce the group homomorphism (weight map) $\wt:\cal P_\mathbb F \to P$ by setting $\wt(\gb\omega_{i,a})=\omega_i$, where $\omega_i$ is the $i$-th fundamental weight of $\lie g$. Otherwise, let $\mathbb K$ be an algebraically closed extension of $\mathbb F$ so that $\cal P_\mathbb F$ can be regarded as a subgroup of $\cal P_\mathbb K$ and define the weight map on $\cal P_\mathbb F$ by restricting the one on $\cal P_\mathbb K$ (this clearly does not depend on the choice of $\mathbb K$).
Define the $\ell$-weight lattice of $U_q(\tlie g)$ to be $\cal P_q:=\cal P_{\mathbb C(q)}$. The submonoid $\cal P_q^+$ of $\cal P_q$ consisting of $n$-tuples of polynomials is called the set of dominant $\ell$-weights of $U_q(\tlie g)$.

Given $\gb\lambda\in\cal P_q^+$ with $\gb\lambda_i(u) = \prod_j (1-a_{i,j}u)$, where $a_{i,j}$ belongs to some algebraic closure of $\mathbb C(q)$, let $\gb\lambda^-\in\cal P_q^+$ be defined by $\gb\lambda^-_i(u) = \prod_j (1-a_{i,j}^{-1}u)$. We will also use the notation $\gb\lambda^+ = \gb\lambda$. Two elements $\gb\lambda,\gb\mu$ of $\cal P_q^+$  are said to be relatively prime if $\gb\lambda_i(u)$ is relatively prime to $\gb\mu_j(u)$ in $\mathbb C(q)[u]$ for all $i,j\in I$. Every  $\gb\nu\in\cal P_q$ can be uniquely written in the form
\begin{equation}\label{e:frac}
\gb\nu = \gb\lambda\gb\mu^{-1} \quad\text{with}\quad \gb\lambda,\gb\mu\in\cal P_q^+ \quad\text{relatively prime}.
\end{equation}
Given $\gb\nu = \gb\lambda\gb\mu^{-1}$ as above, define a $\mathbb C(q)$-algebra homomorphism $\gb\Psi_{\gb\nu}:U_q(\tlie h)\to \mathbb C(q)$ by  setting
\begin{equation}
\gb\Psi_{\gb\nu}(k_i^{\pm 1}) = q_i^{\pm \wt(\gb\nu)(h_i)}, \qquad \sum_{r\ge 0} \gb\Psi_{\gb\nu}(\Lambda_{i,\pm r}) u^r = \frac{(\gb\lambda^{\pm})_i(u)}{(\gb\mu^{\pm})_i(u)}
\end{equation}
where the division is that of formal power series in $u$. The next proposition is easily checked.

\begin{prop}\label{p:lwl*}
The map $\gb\Psi:\cal P_q\to (U_q(\tlie h))^*$ given by $\gb\nu\mapsto \gb\Psi_{\gb\nu}$ is injective.\hfill\qedsymbol
\end{prop}

Define the $\ell$-weight lattice $\cal P$ of $\tlie g$ to be the subgroup of $\cal P_q$ generated by $\gb\omega_{i,a}$ for all $i\in I$ and all $a\in\mathbb C^\times$ or, equivalently, $\cal P = \cal P_\mathbb C$. Observe that every element $\gb\lambda\in \cal P$ can be uniquely decomposed as
\begin{equation}\label{e:fact}
\gb\lambda = \prod_j \gb\omega_{\lambda_j,a_j} \quad\text{for some}\quad \lambda_j\in P \quad\text{and}\quad a_i\ne a_j\in\mathbb C.
\end{equation}
Set also $\cal P^+=\cal P\cap\cal P_q^+$.

From now on we will identify $\cal P_q$  with its image in $(U_q(\tlie h))^*$  under $\gb\Psi$. Similarly, $\cal P$ will be identified with a subset of $U(\tlie h)^*$ via the homomorphism $\gb\Psi_{\gb\nu}:U(\tlie h)\to \mathbb C$ determined by
\begin{equation}
\gb\Psi_{\gb\nu}(h_i) = \wt(\gb\nu)(h_i), \qquad \sum_{r\ge 0} \gb\Psi_{\gb\nu}(\Lambda_{i,\pm r}) u^r = \frac{(\gb\lambda^{\pm})_i(u)}{(\gb\mu^{\pm})_i(u)}.
\end{equation}

It will be convenient to introduce the following notation. Given $i\in I, a\in\mathbb C(q)^\times, r\in\mathbb Z_{\ge 0}$, define
\begin{equation}\label{e:krhlw}
\gb\omega_{i,a,r} = \prod_{j=0}^{r-1} \gb\omega_{i,aq_i^{r-1-2j}}.
\end{equation}
Define also the polynomial
\begin{equation}\label{e:qsf}
f_{i,a,r}(u) = \prod_{j=0}^{r-1} (1-aq_i^{r-1-2j}u).
\end{equation}
Observe that given $f(u)\in\mathbb C(q)[u]$ having all its roots in $\mathbb C(q)$ and such that $f(1)=0$, there exist unique $m\in\mathbb Z_{\ge 0}, a_1,\dots,a_m\in\mathbb C(q)^\times,$ and $r_1,\dots,r_m\in\mathbb Z_{\ge 1}$  such that
\begin{equation}
f(u) = \prod_{k=1}^m f_{i,a_k,r_k}(u)\qquad\text{with}\qquad \frac{a_l}{a_j}\ne q^{\pm(r_l+r_j-2p)}\quad\text{for}\quad 0\le p<\min\{r_l,r_j\}.
\end{equation}
In particular, given $\gb\lambda\in\cal P_q^+$ such that $\gb\lambda_i(u)$ splits in $\mathbb C(q)[u]$ for all $i\in I$, there exist unique $m_i\in\mathbb\mathbb Z_{\ge 0}, a_{i,k}\in\mathbb C(q)^\times$, and $r_{i,k}\in\mathbb Z_{\ge 1}$ such that
\begin{gather}\notag
\gb\lambda = \prod_{i\in I}\prod_{k=1}^{m_i} \gb\omega_{i,a_{i,k},r_{i,k}} \qquad\text{with}\qquad\\ \hfill\\\notag \frac{a_{i,j}}{a_{i,l}}\ne q_i^{\pm(r_{i,j}+r_{i,l}-2p)}\quad\text{and}\quad \sum_{k=1}^{m_i}r_{i,k}=\wt(\gb\lambda)(h_i) \quad\text{for all}\quad i\in I, j\ne l, 0\le p<\min\{r_{i,j},r_{i,l}\}.
\end{gather}

If $J\subseteq I$ and $\gb\lambda\in\cal P_q$, let $\gb\lambda_J$ be the associated $J$-tuple of rational functions. Notice that, if $\gb\lambda_j(u)\in\mathbb C(q_j)(u)$ for all $j\in J$, $\gb\lambda_J$ can be regarded as an element of the $\ell$-weight lattice of $U_q(\tlie g_J)$. Let also $\gb\lambda^J\in\cal P_q$ be such that $(\gb\lambda^J)_j(u)=\gb\lambda_j(u)$ for every $j\in J$ and $(\gb\lambda^J)_j(u)=1$ otherwise.

Recall that $w_0$ defines a Dynkin diagram automorphism such that $w_0\cdot i =j$ iff $w_0\omega_i=-\omega_j$ for  $i,j\in I$.  Given $\gb\lambda\in\cal P_q^+$, let $\gb\lambda^*\in\cal P_q^+$ be the element defined by
\begin{equation}
(\gb\lambda^*)_i(u) = \gb\lambda_{w_0\cdot i}(q^{r^\vee h^\vee}u)
\end{equation}
where $h^\vee$ is the dual Coxeter number of $\lie g$ and  $r^\vee=\max\{c_{ij}c_{ji}:i,j\in I, i\ne j\}$ is the lacing number of $\lie g$. Define also the element $^*\gb\lambda$ by requiring
\begin{equation}
({}^*\gb\lambda)^\pm = (\gb\lambda^*)^{\mp}.
\end{equation}

Given $i\in I$ and $a\in\mathbb C(q)^\times$, define the simple $\ell$-root $\gb\alpha_{i,a}$ by
\begin{equation}
\gb\alpha_{i,a} = (\gb\omega_{i,aq_i,2})^{-1}\prod_{j\ne i} \gb\omega_{j,aq_i,-c_{j,i}}.
\end{equation}
The subgroup of $\cal P_q$ generated by the simple $\ell$-roots is called the $\ell$-root lattice of $U_q(\tlie g)$ and will be denoted by $\cal Q_q$. Let also $\cal Q_q^+$ be the submnoid generated by the simple $\ell$-roots. Quite clearly $\wt(\gb\alpha_{i,a})=\alpha_i$. Define a partial order on $\cal P_q$ by
$$\gb\mu\le\gb\lambda \qquad{if}\qquad \gb\lambda\gb\mu^{-1}\in\cal Q_q^+.$$

\section{Finite-dimensional representations}\label{s:fdimreps}

\subsection{Simple Lie algebras}
We now review some basic facts about the representation theory of $\lie g$ and $U_q(\lie g)$. For the details see \cite{hum:book} and \cite{cp:book} for instance.

Given a $\lie g$-module $V$ and $\mu\in\lie h^*$, let
$$V_\mu=\{v\in V: hv=\mu(h)v \text{ for all } h\in\lie h\}.$$
A nonzero vector $v\in V_\mu$ is called a weight vector of weight $\mu$. If $v$ is a weight vector such that $\lie n^+v = 0$, then $v$ is called a highest-weight vector. If $V$ is generated by a highest-weight vector of weight $\lambda$, then $V$ is said to be a highest-weight module of highest weight $\lambda$.

The following theorem summarizes the basic facts about finite-dimensional $\lie g$-modules.

\begin{thm}\label{t:cig} Let $V$ be a finite-dimensional $\lie g$-module. Then:
\begin{enumerate}
\item $V=\opl_{\mu\in P}^{} V_\mu$ and $\dim V_\mu = \dim V_{w\mu}$ for all $w\in\cal W$.
\item $V$ is completely reducible.
\item For each $\lambda\in P^+$ the $U(\lie g)$-module $V(\lambda)$ generated by a vector $v$ satisfying
$$x_{\alpha_i}^+v=0, \quad h_iv=\lambda(h_i)v, \quad (x_{\alpha_i}^-)^{\lambda(h_i)+1}v=0,\quad\forall\ i\in I,$$
is irreducible and finite-dimensional. If $V$ is irreducible, then
$V$ is isomorphic to $V(\lambda)$ for some $\lambda\in P^+$.
\item If $\lambda\in P^+$ and $V\cong V(\lambda)$, then $V_\mu\ne 0$ iff $w\mu\le\lambda$ for all $w\in\cal W$. Furthermore, the lowest weight of $V(\lambda)$ is $-\lambda^*$. In particular, $V(\lambda)^*\cong V(\lambda^*)$.
\hfill\qedsymbol
\end{enumerate}
\end{thm}

We will need the following lemma.

\begin{lem}\label{l:hwvecs}
Let $V$ be a finite-dimensional $\lie g$-module and suppose $l\in\mathbb Z_{\ge 1}, \nu_k\in P, v_k\in V_{\nu_k}$, for $k=1,\dots,l$, are such that $V=\sum_{k=1}^l U(\lie n^-)v_k$. Fix a decomposition $V= \opl_{j=1}^m V_j$ where $m\in\mathbb Z_{\ge 1},  V_j\cong V(\mu_j)$ for some $\mu_j\in P^+$, and let $\pi_j:V\to V_j$ be the associated projection for $j=1,\dots,m$. Then, there exist distinct $k_1,\dots,k_m\in\{1,\dots,l\}$ such that $\nu_{k_j}=\mu_j$ and $\pi_j(v_{k_j})\ne 0$.
\end{lem}

\begin{proof}
Proceed by induction on $m$. If $m=1$ the lemma is immediate. Otherwise, suppose, without loss of generality, that $\mu_m$ is a maximal weight of $V$. In that case, there must exist $k_m$ such that $\nu_{k_m}=\mu_m$ and $v_{k_m}$ generates an irreducible submodule of $V$ isomorphic to $V(\mu_m)$. In particular, there exists $j$ such that $\mu_j=\mu_m$ and $\pi_j(v_{k_m})\ne 0$. Up to re-ordering, we can assume $j=m$. The lemma now easily follows from the induction hypothesis applied to $\overline V:=V/U(\lie g)v_{k_m}$ and the induced decomposition $\overline V=\opl_{j=1}^{m-1} \overline V_j$ where $\overline V_j$ is the image of $V_j$ in $\overline V$.
\end{proof}

Let $\mathbb Z[P]$ be the integral group ring over $P$ and denote by $e:P\to\mathbb Z[P], \lambda\mapsto e^\lambda$, the inclusion of $P$ in $\mathbb Z[P]$ so that $e^\lambda e^\mu=e^{\lambda+\mu}$. Given a finite-dimensional $\lie g$-module $V$, the character of $V$ is defined to be
\begin{equation}\label{e:chd}
\ch(V) = \sum_{\mu\in P} \dim(V_\mu) e^\mu.
\end{equation}

Given a $U_q(\lie g)$-module $V$ and $\mu\in P$, let
$$V_\mu=\{v\in V: k_iv=q_i^{\mu(h_i)}v \text{ for all } i\in I\}.$$
A nonzero vector $v\in V_\mu$ is called a weight vector of weight $\mu$. If $v$ is a weight vector such that $x_i^+v = 0$ for all $i\in I$, then $v$ is called a highest-weight vector. If $V$ is generated by a highest-weight vector of weight $\lambda$, then $V$ is said to be a highest-weight module of highest weight $\lambda$.
A $U_q(\lie g)$-module $V$ is said to be a weight module if $V=\opl_{\mu\in P}^{} V_\mu$.  Denote by $\cal C_q$ be the category of all finite-dimensional weight modules of $U_q(\lie g)$.

\begin{rem}
A $U_q(\lie g)$-module $V$ satisfying $V=\opl_{\mu\in P}^{} V_\mu$ is usually called a weight-module of type {\bf 1}. We shall not discuss what type {\bf 1} means here. For further details see \cite{cp:book} for instance.
\end{rem}

The character of an object $V\in\cal C_q$ is defined by \eqref{e:chd}.
The following theorem is the quantum analogue of Theorem \ref{t:cig}.

\begin{thm}\label{t:ciuqg} Let $V\in\cal C_q$. Then:
\begin{enumerate}
\item $\dim V_\mu = \dim V_{w\mu}$ for all $w\in\cal W$.
\item $V$ is completely reducible.
\item For each $\lambda\in P^+$ the $U(\lie g)$-module $V_q(\lambda)$ generated by a vector $v$ satisfying
$$x_i^+v=0, \quad k_iv=q_i^{\lambda(h_i)}v, \quad (x_i^-)^{\lambda(h_i)+1}v=0,\quad\forall\ i\in I,$$
is irreducible and finite-dimensional. If $V$ is irreducible, then
$V$ is isomorphic to $V_q(\lambda)$ for some $\lambda\in P^+$.
\item If $\lambda\in P^+$ and $V\cong V_q(\lambda)$, then $\ch(V) = \ch(V(\lambda))$. In particular, $V_q(\lambda)^*\cong V_q(\lambda^*)$.
\hfill\qedsymbol
\end{enumerate}
\end{thm}

If $J\subseteq I$ we shall denote by $V_q(\lambda_J)$ the $U_q(\lie g_J)$-irreducible module of highest weight $\lambda_J$. Similarly $V(\lambda_J)$ denotes the corresponding irreducible $\lie g_J$-module

\subsection{Loop algebras}

Let $V$ be a $U_q(\tlie g)$-module. We say that a nonzero vector $v\in V$ is an $\ell$-weight vector if there exists $\gb\lambda\in\cal P_q$ and $k\in\mathbb Z_{>0}$ such that $(\eta-\gb\Psi_\gb\lambda(\eta))^kv=0$ for all $\eta\in U_q(\tlie h)$. In that case, $\gb\lambda$ is said to be the $\ell$-weight of $v$. $V$ is said to be an $\ell$-weight module if every vector of $V$ is a linear combination of $\ell$-weight vectors. In that case, let $V_\gb\lambda$ denote the subspace spanned by all $\ell$-weight vectors of $\ell$-weight $\gb\lambda$.
An $\ell$-weight vector $v$ is said to be a highest-$\ell$-weight vector if $\eta v=\gb\Psi_\gb\lambda(\eta)v$ for every $\eta\in U_q(\tlie h)$ and $x_{i,r}^+v=0$ for all $i\in I$ and all $r\in\mathbb Z$. $V$ is said to be a highest-$\ell$-weight module if it is generated by a highest-$\ell$-weight vector. The notion of lowest-$\ell$-weight module is defined similarly. Denote by $\wcal C_q$ the category of all finite-dimensional $\ell$-weight modules of $U_q(\tlie g)$.
Quite clearly $\wcal C_q$ is an abelian category.

Observe that if $V\in\wcal C_q$, then $V\in\cal C_q$ and
\begin{equation}\label{e:lws}
V_\lambda = \opl_{\gb\lambda:\wt(\gb\lambda)=\lambda}^{} V_\gb\lambda.
\end{equation}
Moreover, if $V$ is a highest-$\ell$-weight module of highest $\ell$-weight $\gb\lambda$, then
\begin{equation}\label{e:hw}
\dim(V_{\wt(\gb\lambda)}) = 1\qquad\text{and}\qquad V_\mu\ne 0 \Rightarrow \mu\le\wt(\gb\lambda).
\end{equation}

Define the concepts of $\ell$-weight vector, etc., for $\tlie g$ in a similar way  and denote by $\wcal C$ the category of all finite-dimensional $\tlie g$-modules. The next proposition is easily established using \eqref{e:hw}.

\begin{prop}
If $V$ is a highest-$\ell$-weight module, then it has a unique proper submodule and, hence, a unique irreducible quotient.\hfill\qedsymbol
\end{prop}

\begin{defn}
Let $\gb\lambda\in \cal P_q^+$ and $\lambda=\wt(\gb\lambda)$. The Weyl module $W_q(\gb\lambda)$ of highest $\ell$-weight $\gb\lambda$ is the $U_q(\tlie g)$-module defined by the quotient of $U_q(\tlie g)$ by the left ideal generated by the elements $x_{i,r}^+, (x_{i,r}^-)^{\lambda(h_i)+1}$, and $\eta-\gb\Psi_\gb\lambda(\eta)$ for every $i\in I, r\in\mathbb Z$, and $\eta\in U_q(\tlie h)$. Denote by $V_q(\gb\lambda)$ the irreducible quotient of $W_q(\gb\lambda)$. The Weyl module $W(\gb\lambda), \gb\lambda\in \cal P^+$, of $\tlie g$ is defined in a similar way. Its irreducible quotient will be denoted by $V(\gb\lambda)$.
\end{defn}

The next theorem was proved in \cite[Lemmas 4.6 and 4.7]{cp:weyl} for simply laced $\lie g$ and in \cite[Proposition 2.2]{cha:fer} for $\lie g$ with lacing number $r^\vee=2$. For the sake of completeness, a proof for $\lie g$ of type $G_2$ will appear in \cite{jm:root}.

\begin{thm}\label{t:weyl}
For every $\gb\lambda\in\cal P_q^+$ (resp. $\cal P^+$) the module $W_q(\gb\lambda)$ (resp. $W(\gb\lambda)$) is the universal finite-dimensional $U_q(\tlie g)$-module (resp. $\tlie g$-module) with highest $\ell$-weight $\gb\lambda$. Every simple object of $\wcal C_q$ (resp. $\wcal C$) is highest-$\ell$-weight. \hfill\qedsymbol
\end{thm}

\begin{rem}
It is not true that the module $V_q(\gb\lambda)$ belongs to $\wcal C_q$ for every $\gb\lambda\in\cal P_q^+$. This is so because $\mathbb C(q)$ is not algebraically closed. In fact, one can prove, using some results of subsection \ref{ss:qchar} below, that $V_q(\gb\lambda)$ is in $\wcal C_q$ iff $\gb\lambda_i(u)$ splits in $\mathbb C(q)[u]$ for every $i\in I$. Otherwise, $V_q(\gb\lambda)$ is quasi-$\ell$-weight in a sense analogous to that defined in \cite{jm:nacf} in the context of hyper loop algebras.
\end{rem}

We shall need the following lemma which is a consequence of the proof of Theorem \ref{t:weyl}.

\begin{lem}\label{l:curalggen}
If $V$ is a highest-$\ell$-weight module of $\tlie g$ and $v$ be a highest-$\ell$-weight vector. Then $V=U(\lie g[t])v$.
\end{lem}

If $J\subseteq I$ we shall denote by $V_q(\gb\lambda_J)$ the $U_q(\tlie g_J)$-irreducible module of highest $\ell$-weight $\gb\lambda_J$. Similarly $V(\gb\lambda_J)$ denotes  the corresponding irreducible $\tlie g_J$-module. Similar notations for the Weyl modules are defined in the obvious way.

We shall need the following result about dual representations proved in \cite{freem:qchar}.

\begin{prop}\label{p:dual}
For every $\gb\lambda\in\cal P_q^+$, $V_q(\gb\lambda)$ is a lowest-$\ell$-weight module with lowest $\ell$-weight $(\gb\lambda^*)^{-1}$. In particular,$V_q(\gb\lambda)^* \cong V_q(\gb\lambda^*)$.\hfill\qedsymbol
\end{prop}

\subsection{Evaluation modules and Cartan involution}

Given a $\lie g$-module $V$, let $V(a)$ be the $\tlie g$-module obtained by pulling-back the evaluation map ${\rm ev}_a$. Such modules are called evaluation modules. If $V=V(\lambda)$ we use the notation $V(\lambda,a)$ for the corresponding evaluation module. The next theorem was proved in \cite{cha:int,cp:new,cp:weyl}.

\begin{thm}\label{t:tev}
Let $\gb\lambda\in\cal P^+$.
\begin{enumerate}
\item If $\gb\lambda = \gb\omega_{\lambda,a}$ for some $\lambda\in P^+$ and some $a\in\mathbb C^\times$, then $V(\gb\lambda)\cong V(\lambda,a)$.
\item If $\gb\lambda = \prod_j \gb\omega_{\lambda_j,a_j}$ as in \eqref{e:fact}, then $V(\gb\lambda)\cong \otm_j^{} V(\lambda_j,a_j)$ and $W(\gb\lambda)\cong\otm_j^{} W(\gb\omega_{\lambda_j,a_j})$.\hfill\qedsymbol
\end{enumerate}
\end{thm}

\begin{cor}
Every object in $\wcal C$ is an $\ell$-weight module.\hfill\qedsymbol
\end{cor}

Assume $\lie g$ is of type $A$ and consider the $\mathbb C(q)$-algebra $U_q'(\lie g)$ given by generators $x_i^\pm, k_\mu^{\pm 1}$ with $i\in I, \mu\in P$, and the following defining relations:
\begin{gather*}
k_\mu k_\mu^{-1} = k_\mu^{-1}k_\mu =1, \qquad k_\mu k_\nu = k_{\mu+\nu}\\
k_\mu x_{j}^\pm k_\mu^{-1} = q^{\mu(h_j)}x_{j}^{\pm},\qquad  [x_{i}^+ , x_{j}^-]=\delta_{i,j} \frac{k_{\alpha_i} - k_{\alpha_i}^{-1}}{q - q^{-1}},\\
\sum_{k=0}^{1-c_{ij}}(-1)^k   (x_i^\pm)^{(1-c_{ij}-k)} x_{j}^{\pm} (x_i^\pm)^{(k)} =0,\ \ \text{if $i\ne j$},
\end{gather*}
There is an obvious monomorphism of algebras $U_q(\lie g)\to U_q'(\lie g)$ such that $k_i\mapsto k_{\alpha_i}$.
A $U_q'(\lie g)$-module is said to be a weight module if the generators $k_\nu, \nu\in P$, act diagonally with eigenvalues of the form $q^{(\nu,\mu)}$ for some $\mu\in P$ where $(\cdot,\cdot)$ is the bilinear form such that $(\alpha_i,\alpha_j)=c_{ij}$. It is not difficult to to see that restriction establishes an equivalence of categories from that $U_q'(\lie g)$-weight modules to $\cal C_q$. From now on we shall identify these two categories using this equivalence. The next proposition was proved in \cite[\S2]{jim:qan} and \cite[Proposition 3.4]{cp:small}.

\begin{prop}\label{p:qev}
Let $\lie g$ be of type $A$. Then, there exists an algebra homomorphism ${\rm qev}: U_q(\tlie g)\to U_q'(\lie g)$ satisfying: if $\lambda\in P^+$ and $V$ is the pull-back of $V_q(\lambda)$ by ${\rm qev}$, then there exists $l(\lambda)\in\mathbb Z$ such that $V$ is isomorphic to $V_q(\gb\lambda)$ where
$$\gb\lambda = \prod_{i\in I}\gb\omega_{i,a_i,\lambda(h_i)} \qquad\text{with}\qquad a_1 = q^{l(\lambda)} \quad\text{and}\quad \frac{a_{i+1}}{a_i} = q^{\lambda(h_i)+\lambda(h_{i+1})+1} \quad\text{for}\quad i<n.$$

\vspace{-33pt}\hfill\qedsymbol
\end{prop}
\vspace{10pt}
Given  $a\in\mathbb C(q)^\times$, there exists a unique $\mathbb C(q)$-algebra automorphism $\varrho_a$ of $U_q(\tlie g)$ such that $\varrho_a$ is the identity on $U_q(\lie g)$ and $\varrho_a(x_{i,r}^\pm) = a^r x_{i,r}^\pm$. Let
\begin{equation}
{\rm qev}_a = {\rm qev}\circ\varrho_a.
\end{equation}
Denote by $V_q(\lambda,a)$ the pull-back of $V_q(\lambda)$ by the evaluation map ${\rm qev}_a$. It is easy to see that $V_q(\lambda,a)\cong V_q(\gb\lambda) $ where
\begin{equation*}
\gb\lambda = \prod_{i\in I}\gb\omega_{i,a_i,\lambda(h_i)} \qquad\text{with}\qquad a_1 = aq^{l(\lambda)} \quad\text{and}\quad \frac{a_{i+1}}{a_i} = q^{\lambda(h_i)+\lambda(h_{i+1})+1} \quad\text{for}\quad i<n.
\end{equation*}

It turns out that, for $\lie g$ not of type $A$, there is no analogue of the map ${\rm qev}$. In fact, it is known (see \cite{cha:fer} for instance) that there exists $i\in I$ and $m\in\mathbb Z_{\ge 0}$ such that the action of $U_q(\lie g)$ on $V_q(m\omega_i)$ cannot be extended to one of $U_q(\tlie g)$.

One easily checks that there exists a unique algebra involution $\tilde\sigma$ of $U_q(\tlie g)$ such that $\tilde\sigma(x_{i,r}^\pm) = x_{i,-r}^{\mp}, \tilde\sigma(k_i) = k_i^{-1}$, and $\tilde\sigma(h_{i,s})=-h_{i,-s}$ for all $i\in I, r,s\in\mathbb Z, s\ne 0$. The involution $\tilde\sigma$ is called Cartan involution and it is also a coalgebra anti-involution. The restriction of $\tilde\sigma$ to $U_q(\lie g)$ defines an involution $\sigma$ of $U_q(\lie g)$ also called Cartan involution. Given a $U_q(\tlie g)$-module $V$, let $V^{\tilde\sigma}$ be the pullback of $V$ by $\tilde\sigma$. Similarly, $V^\sigma$ will denote the pullback of a $U_q(\lie g)$-module $V$ by $\sigma$. It is not difficult to see that a highest-$\ell$-weight vector of $V_q(\gb\lambda)$ is a lowest-$\ell$-weight vector of $V_q(\gb\lambda)^{\tilde\sigma}$. Moreover, it follows from \eqref{e:Lambdad} that
\begin{equation}
\tilde\sigma(\Lambda_i^\pm(u)) = (\Lambda_i^{\mp}(u))^{-1}
\end{equation}
where
\begin{equation}
\Lambda_i^\pm(u)=\sum_{r=0}^\infty \Lambda_{i,\pm r} u^r
\end{equation}
and the inverse is that of formal power series in $u$. It is now not difficult to complete the proof of the next proposition.

\begin{prop}
Let $\lambda\in P^+$ and $\gb\lambda\in\cal P_q^+$. Then, $V_q(\lambda)^\sigma\cong V_q(\lambda^*)$ and $V_q(\gb\lambda)^{\tilde\sigma}\cong V_q({}^*\gb\lambda)$.\hfill\qedsymbol
\end{prop}
The analogous result in the classical case is established similarly.

We end this subsection by remarking the following. Let $\lie g$ be of type $A$, suppose $\gb\lambda\in\cal P_q^+$ is such that $V_q(\gb\lambda)\cong V_q(\lambda,a)$ for some $a\in\mathbb C(q)^\times$,  and set $b_{n} = (aq^{l(\lambda)+n+1})^{-1}$. Then,
\begin{equation}\label{e:*lambdaA}
^*\gb\lambda = \prod_{i\in I}\gb\omega_{i,b_i,\lambda^*(h_i)} \qquad\text{with}\qquad  \frac{b_i}{b_{{i-1}}} = q^{-(\lambda^*(h_i)+\lambda^*(h_{i-1})+1)} \quad\text{for}\quad i>1.
\end{equation}

\subsection{Classical limits}

\begin{defn}
Denote by $\cal P_\mathbb A^+$ the subset of $\cal P_q$ consisting of $n$-tuples of polynomials with coefficients in $\mathbb A$. Let also $\cal P_\mathbb A^{++}$ be the subset of $\cal P_\mathbb A^+$ consisting of $n$-tuples of polynomials whose leading terms are in $\mathbb Cq^{\mathbb Z}\backslash\{0\}=\mathbb A^\times$. Given $\gb\lambda\in\cal P_\mathbb A^+$, let \bgb\lambda\ be the element of $\cal P^+$ obtained from $\gb\lambda$ by evaluating $q$ at $1$.
\end{defn}

Recall that an $\mathbb A$-lattice (or form) of a $\mathbb C(q)$-vector space $V$ is a free $\mathbb A$-submodule $L$ of $V$ such that $\mathbb C(q)\otimes_\mathbb A L=V$. If $V$ is a $U_q(\tlie g)$-module, a $U_\mathbb A(\tlie g)$-admissible lattice of $V$ is an $\mathbb A$-lattice of $V$ which is also a $U_\mathbb A(\tlie g)$-submodule of $V$. Given a $U_\mathbb A(\tlie g)$-admissible lattice of a $U_q(\tlie g)$-module $V$, define
\begin{equation}\label{e:clm}
\bar L = \mathbb C\otimes_\mathbb A L,
\end{equation}
where $\mathbb C$ is regarded as an $\mathbb A$-module by letting $q$ act as $1$. Then $\bar L$ is a $\tlie g$-module by Proposition \ref{p:cluq} and $\dim(\bar L)=\dim(V)$. The next theorem is essentially a corollary of the proof of Theorem \ref{t:weyl}.

\begin{thm}\label{t:llattice}
Let $V$ be a nontrivial quotient of $W_q(\gb\lambda)$ for some $\gb\lambda\in \cal P_\mathbb A^{++}$, $v$  a highest-$\ell$-weight vector of $V$, and $L=U_\mathbb A(\tlie g)v$. Then, $L$ is a $U_\mathbb A(\tlie g)$-admissible lattice of $V$ and $\ch(\bar L)=\ch(V)$. In particular, $\bar L$ is a quotient of $W(\bgb\lambda)$.\hfill\qedsymbol
\end{thm}

\begin{defn}
Let $\gb\lambda\in \cal P_\mathbb A^{++}$, $v$ be a highest-$\ell$-weight vector of $V_q(\gb\lambda)$ and $L=U_\mathbb A(\tlie g)v$. We denote by $\overline{V_q(\gb\lambda)}$ the $\tlie g$-module $\bar L$.
\end{defn}

We shall also use the following straightforward lemma

\begin{lem}\label{l:morlim}
Let $V,V'$ be $U_q(\tlie g)$-modules and $L,L'$ be $U_\mathbb A(\tlie g)$-submodules. Suppose $\phi:V\to V'$ is a $U_q(\tlie g)$-module map such that $\phi(L)\subseteq L'$. Then $\bar\phi:=1\otimes \phi : \bar L\to\bar L'$ is a $\tlie g$-module map.\hfill\qedsymbol
\end{lem}

\section{Minimal affinizations}\label{s:min}

\subsection{Classification} We now review the notion of minimal affinizations of an irreducible $U_q(\lie g)$-module introduced in \cite{cha:minr2}.

Given $\lambda\in P^+$, an object $V\in\wcal C_q$ is said to be an affinization of $V_q(\lambda)$ if, as a $U_q(\lie g)$-module,
\begin{equation}
V\cong V_q(\lambda)\oplus \opl_{\mu< \lambda}^{} V_q(\mu)^{\oplus m_\mu(V)}
\end{equation}
for some $m_\mu(V)\in\mathbb Z_{\ge 0}$. Two affinizations of $V_q(\lambda)$ are said to be equivalent if they are isomorphic as $U_q(\lie g)$-modules. If $\gb\lambda\in\cal P_q^+$ is such that $\wt(\gb\lambda)=\lambda$, then $V_q(\gb\lambda)$ is quite clearly an affinization of $V_q(\lambda)$. The partial order on $P^+$ induces a natural partial order on the set of  (equivalence classes of) affinizations of $V_q(\lambda)$. Namely, if $V$ and $W$ are affinizations of $V_q(\lambda)$, say that $V\le W$ if either $m_\mu(V)\le m_\mu(W)$ for all $\mu\in P^+$ or if for all $\mu\in P^+$ such that $m_\mu(V)>m_\mu(W)$ there exists $\nu>\mu$ such that $m_\nu(V)<m_\nu(W)$.
A minimal element of this partial order is said to be a minimal affinization.

Suppose $\lie g$ is not of types $D$ or $E$. Given $\gb\lambda\in\cal P_q^+$ set
\begin{equation}
\gb\lambda^o = \gb\lambda^* \quad \text{if}\quad\lie g=\lie{sl}_{n+1} \qquad\text{and}\qquad \gb\lambda^o = {}^*\gb\lambda \quad\text{otherwise.}
\end{equation}
Recall that, in these cases, $\lambda^*=\lambda$ for all $\lambda\in P^+$ except if $\lie g$ is of type $A$.

The following is the main result of \cite{cha:minr2,cp:minsl,cp:minnsl} and it gives a partial classification of the highest $\ell$-weights of the minimal affinizations. In fact it gives the complete classification when $\lie g$ is not of types $D$ or $E$.

\begin{thm}\label{t:drimin}
Let $\gb\lambda\in\cal P^+_q,\lambda=\wt(\gb\lambda)$, and  $V=V_q(\gb\lambda)$. Suppose $\lie g$ is not of types $D$ or $E$. Then $V$ is a minimal affinization of $V_q(\lambda)$ iff $V^*$ and $V^{\tilde\sigma}$ are minimal affinizations of $V_q(\lambda^*)$. In that case, there exist $a\in\mathbb C(q)^\times$ and $\mu\in \{\lambda,\lambda^*\}$ such that  either $\gb\lambda$ or $\gb\lambda^o$ is equal to
$$\prod_{i=1}^n \gb\omega_{i,a_i,\mu(h_i)} \qquad\text{with}\qquad a_1 = a \qquad\text{and}\qquad \frac{a_{i+1}}{a_i} = q^{d_i\mu(h_i)+d_{i+1}\mu(h_{i+1})+r_{i}^\vee}$$
for all $i\in I, i<n$, where $r_i^\vee = d_{i}-1-c_{i,i+1}$. Equivalently, $V$ is a minimal affinization of $V_q(\lambda)$ iff there exist $a\in\mathbb C(q)^\times$ and $\epsilon\in\{1,-1\}$ such that
$$\gb\lambda=\prod_{i=1}^n \gb\omega_{i,a_i,\lambda(h_i)} \qquad\text{with}\qquad a_1 = a \qquad\text{and}\qquad \frac{a_{i+1}}{a_i} = q^{\epsilon(d_i\lambda(h_i)+d_{i+1}\lambda(h_{i+1})+r_{i}^\vee)}$$
for all $i\in I, i<n$.
If $\lie g$ is of type $D$ or $E$, suppose the support of $\lambda$ is contained in a connected subdiagram $J\subseteq I$ of type $A$. Then, $V$ is a minimal affinization of $V_q(\lambda)$ iff $V_q(\gb\lambda_J)$ is a minimal affinization of $V_q(\lambda_J)$. \hfill\qedsymbol
\end{thm}

The next corollary is immediate (recall from \S\ref{ss:clalg} that $\overline{\rm supp}(\lambda)$ is the minimal connected subdiagram of $I$ containing ${\rm supp}(\lambda)$).

\begin{cor}
Suppose $\lambda\in P^+$ is such that $\overline{\rm supp}(\lambda)$ does not contain a subdiagram of type $D_4$. Then, $V_q(\lambda)$ has a unique equivalence class of minimal affinizations.\hfill\qedsymbol
\end{cor}

\begin{rem}
We warn the reader that the conditions we give in Theorem \ref{t:drimin} do not match the ones given in \cite{cha:minr2,cp:minsl,cp:minnsl}. This is due to different normalizations in some definitions. Our notation follows more closely that of \cite{her:min} which is more uniform. We also notice that $r_i^\vee=d_{i+1}-c_{i+1,i}$ and, moreover, $r_i^\vee\in\{r^\vee-1,r^\vee\}$ for all $i\in I, i<n$. It is easy to check that $r_i^\vee=r^\vee$ for all $i<n$ if $\lie g$ is of types $A, B$, or $G$. If $\lie g$ is of type $C$, then $r_i^\vee=r^\vee-1$ iff $i<n-1$. Finally, if $\lie g$ is of type $F$, then $r_i^\vee=r^\vee$ iff $\alpha_i$ is a long root.
\end{rem}

\begin{cor}
For every $a\in\mathbb C(q)^\times, i\in I$ and $m\in\mathbb Z_{\ge 0}$, the module $V_q(\gb\omega_{i,a,m})$ is a minimal affinization of $V_q(m\omega_i)$.\hfill\qedsymbol
\end{cor}

The modules $V_q(\gb\omega_{i,a,m})$ are known as Kirillov-Reshetikhin modules.

In the cases not covered by Theorem \ref{t:drimin}, i.e., when $\overline{\rm supp}(\lambda)$ contains a subdiagram of type $D_4$, then $V_q(\lambda)$ may have more then one equivalence class of minimal affinizations (see \cite{cp:minsl,cp:minirr}). We shall briefly discuss these cases in sections \ref{ss:multmar} and \ref{ss:multmai}.

We now state a few results which were used in the proof of Theorem \ref{t:drimin} and will be useful for us as well. The proofs can be found in \cite{cp:minsl}.

\begin{lem}\label{l:dsa}
Suppose $\emptyset\ne J\subseteq I$ is a connected subdiagram of the Dynkin diagram of $\lie g$. Let $V=V_q(\gb\lambda)$, $v$ a highest-$\ell$-weight vector of $V$, and $V_J = U_q(\tlie g_J)v$. Then, $V_J\cong V_q(\gb\lambda_J)$.\hfill\qedsymbol
\end{lem}

\begin{defn}
Suppose $\lie g$ is not of type $D$ or $E$. A connected subdiagram $J\subseteq I$ is said to be admissible if $J$ is of type $A$. If $\lie g$ is of type $D$ or $E$, let $i_0\in I$ be the unique element connected to three other nodes. A connected subdiagram $J\subseteq I$ is said to be admissible if $J$ is of type $A$ and $J\backslash\{i_0\}$ is connected.
\end{defn}

\begin{prop}\label{p:admsd}
Suppose $J\subseteq I$ is admissible and that $\gb\lambda\in\cal P_q^+$ is such that $V_q(\gb\lambda)$ is a minimal affinization of $V_q(\lambda)$ where $\lambda=\wt(\gb\lambda)$. Then $V_q(\gb\lambda_J)$ is a minimal affinization of $V_q(\lambda_J)$.\hfill\qedsymbol
\end{prop}

\begin{prop}\label{p:fracqz}
Let $\gb\lambda\in P_q^+$ and $\lambda=\wt(\gb\lambda)$. If $V_q(\gb\lambda)$ is a minimal affinization of $V_q(\lambda)$, then there exist $a_i\in\mathbb C(q)^\times,i\in I$, such that $\gb\lambda = \prod_{i\in I}\gb\omega_{i,a_i,\lambda(h_i)}$ and $\frac{a_i}{a_j}\in q^{\mathbb Z}$ for all $i,j\in I$.
\end{prop}

\begin{proof}
The existence of $a_i\in\mathbb C(q)^\times,i\in I$, such that $\gb\lambda = \prod_{i\in I}\gb\omega_{i,a_i,\lambda(h_i)}$ follows from Proposition \ref{p:admsd} together with the classification of minimal affinizations for $\lie{sl}_2$. The condition $\frac{a_i}{a_j}\in q^{\mathbb Z}$ for all $i,j\in I$, can be proved using the results of \cite{cha:braid} (cf. subsection \ref{ss:tpkr} below).  Alternatively, the proposition is immediate from Theorems  \ref{t:drimin} and \ref{t:regmma} in the cases covered by them.
\end{proof}

\begin{cor}
For every $\lambda\in P^+$ there exist $\gb\lambda\in\cal P_\mathbb A^{++}$ such that $V_q(\gb\lambda)$ is a minimal affinization of $V_q(\lambda)$.\hfill\qedsymbol
\end{cor}

\subsection{Restricted limits}
In this subsection we define ``restricted limits'' of minimal affinizations. We begin with the following definition.

\begin{defn}
Let $V$ be a $\mathbb Z_{\ge 0}$-graded vector space and denote its $s$-th graded piece by $V[s]$. A $\lie g[t]$-module $V$ is said to be $\mathbb Z_{\ge 0}$-graded if $V$ is a $\mathbb Z_{\ge 0}$-graded vector space and $x\otimes t^rv\in V[r+s]$ for every $v\in V[s], x\in\lie g, r,s\in\mathbb Z_{\ge 0}$. A $\mathbb Z_{\ge 0}$-graded $\lie g[t]$-module $V$ satisfying $V[r]=0$ for $r\gg 0$ is said to be a restricted $\lie g[t]$-module. If $V$ is a $\mathbb Z_{\ge 0}$-graded $\lie g[t]$-module, denote by $V(s)$ the quotient of $V$ by its $\lie g[t]$-submodule $\opl_{k> s}^{} V[k]$.
\end{defn}

The next lemma follows immediately from Proposition \ref{p:fracqz}.

\begin{lem}\label{l:limhlw}
Suppose $\gb\lambda\in\cal P_\mathbb A^{++}$ is such that $V_q(\gb\lambda)$ is a minimal affinization. Then $\bgb\lambda = \gb\omega_{\lambda,a}$ for some $a\in\mathbb C^\times$, where $\lambda=\wt(\gb\lambda)$. \hfill\qedsymbol
\end{lem}

\begin{prop}\label{p:admevrel}
Suppose $\gb\lambda\in\cal P_\mathbb A^{++}$ is such that $V_q(\gb\lambda)$ is a minimal affinization and that $J\subseteq I$ is an admissible subdiagram. Let $v$ be a highest-$\ell$-weight vector of $V=\overline{V_q(\gb\lambda)}, \lambda=\wt(\gb\lambda)$, and $a\in\mathbb C^\times$ be such that $\bgb\lambda = \omega_{\lambda,a}$. Then $x_{\alpha,r}^-v=a^rx_\alpha^-v$ for every $\alpha\in R^+_J$.
\end{prop}

\begin{proof}
Let $J$ be admissible, $\alpha\in R^+_J$, and $V_J=U_q(\tlie g_J)v'$ where $v'\in V_q(\gb\lambda)$ is such that $\overline{v'}=v$. Then $V_J$ is a minimal affinization by Proposition \ref{p:admsd} and, since $J$ is of type $A$, $V_J$ is irreducible as a $U_q(\lie g_J)$-module by Theorem \ref{t:drimin}. Hence, the $\tlie g_J$-submodule of $V$ generated by $v$ is isomorphic to $V(\lambda_J,a)$.
\end{proof}

Recall the definition of the maps $\tau_a:\lie g[t]\to\lie g[t]$ from subsection \ref{ss:clalg}.

\begin{defn}
Let $\gb\lambda\in\cal P_\mathbb A^{++}, \lambda=\wt(\gb\lambda)$, and $a\in\mathbb C^\times$ be such that $\bgb\lambda = \omega_{\lambda,a}$. The $\lie g[t]$-module $L(\gb\lambda)$ is defined to be the pullback of $\overline{V_q(\gb\lambda)}$ by $\tau_a$.  Define also the module $A(\lambda)$ to be the $\lie g[t]$-module given by the quotient of $U(\lie g[t])$ by the left ideal generated by
\begin{equation*}
\lie n^+[t], \qquad \lie h\otimes t\mathbb C[t], \qquad h_i-\lambda(h_i), \qquad (x_{\alpha_i}^-)^{\lambda(h_i)+1}, \qquad x_{\alpha,1}^-
\end{equation*}
for all $i\in I$ and all $\alpha\in R^+_J$ for some admissible subdiagram $J\subseteq I$. Denote by $v_\lambda$ the image of $1$ in $A(\lambda)$ so that $A(\lambda)=U(\lie n^-[t])v_\lambda$.
\end{defn}

It immediately follows from Theorem \ref{t:tev}, Proposition \ref{p:admevrel}, and Lemma \ref{l:curalggen} that $L(\gb\lambda)$ is a quotient of $A(\lambda)$. It is also clear that $A(\lambda)$ is a $\mathbb Z_{\ge 0}$-graded $\lie g[t]$-module. We call the module $L(\gb\lambda)$ the restricted limit of $V_q(\gb\lambda)$. It is immediate from Theorem \ref{t:llattice} that
\begin{equation}\label{e:Misoclass}
\ch(L(\gb\lambda))=\ch(V_q(\gb\lambda)).
\end{equation}
In the special case that $\lambda=m\omega_i$ for some $m\in\mathbb Z_{\ge 0}$ and some $i\in I$, the modules $L(\gb\omega_{i,a,m})$ are called the restricted Kirillov-Reshetikhin modules of highest-weight $m\omega_i$. For $\lie g$ of classical type they were studied in \cite{cm:kr} and for $\lie g$ of type $G_2$ they were studied in \cite{cm:krg}.

\begin{prop}\label{p:Afdim}
For every $\lambda\in P^+$ the module $A(\lambda)$ is finite-dimensional. In particular, $A(\lambda)$ is restricted.
\end{prop}

\begin{proof}
Since $A(\lambda)=U(\lie n^-[t])v_\lambda$. It immediately follows that $(A(\lambda)[r])_\mu$ is finite-dimensional for every $r\in\mathbb Z_{\ge 0}$ and every $\mu\in P$. The relations $(x_{\alpha_i}^-)^{\lambda(h_i)+1}v_\lambda=0$ for all $i\in I$ implies, as usual, that the elements $x_{\alpha_i}^\pm$ act locally nilpotently on $A(\lambda)$ and, hence, $\dim(A(\lambda)_\mu)=\dim(A(\lambda)_{w\mu})$ for every $\mu\in P$ and $w\in \cal W$. This in turn implies that $A(\lambda)_\mu\ne 0$ iff $w_0\lambda\le\mu\le\lambda$. Hence, $A(\lambda)$ has only finitely-many non-trivial weight spaces.
Using the defining relations of $A(\lambda)$ together with basic commutation relations in $\lie g[t]$, it is trivial to see that $x_{\alpha,r}^-v_\lambda=0$ for all $\alpha\in R^+$ provided $r\gg 0$. This together with the PBW theorem then implies that $(A(\lambda)[s])_\mu=0$ for every $\mu\in P$ provided $s\gg 0$. In fact, let $r\in\mathbb Z_{\ge 0}$ be such that $x_{\alpha,s}^-v_\lambda=0$ for all $\alpha\in R^+$ and all $s\ge r$. Fix a total order on $R^+\times\mathbb Z_{\ge 0}$ such that $(\alpha,k)<(\beta,l)$ whenever $k<r$ and $l\ge r$. The PBW monomials for $U(\lie n^-[t])$ are then formed such that $x_{\beta,l}^-$ occur to the right of  $x_{\alpha,k}^-$ whenever $(\alpha,k)<(\beta,l)$. Hence, in order to get to the $s$-th graded piece of $A(\lambda)$ with $s\gg r$, one would have to apply elements of the form $x_{\alpha,k}^-$ with $k<r$  to $v_\lambda$ ``too many times''. This implies that the maximal possible weight of $A(\lambda)[s]$ would fall out of the set of weights lying in between $w_0\lambda$ and $\lambda$.
\end{proof}

\subsection{Relations for $L(\gb\lambda)$}

We now state our main results and conjectures.

\begin{defn}\label{d:KRres}
Let $m\in \mathbb Z_{\ge 0}$ and $i\in I$. The $\lie g[t]$-module $M(m\omega_i)$ is the quotient of $U(\lie g[t])$ by the left ideal generated by
\begin{equation}\label{e:NKRdef}
\lie n^+[t], \qquad \lie h\otimes t\mathbb C[t], \qquad h_j,\qquad h_i-m,\qquad x_{\alpha_j}^-, \qquad (x_{\alpha_i}^-)^{m+1},\qquad x_{\alpha_i,1}^-\qquad \text{for all } j\ne i.
\end{equation}
\end{defn}

Quite clearly $M(m\omega_i)$ is a $\mathbb Z_{\ge 0}$-graded $\lie g[t]$-module and $A(m\omega_i)$ is a quotient of $M(m\omega_i)$. The next proposition follows from the results of \cite{cha:fer,cm:kr,cm:krg}.

\begin{prop}\label{p:cmkr}
Suppose $\lie g$ is not of type $E$ or $F$. Let $i\in I, m\in\mathbb Z_{\ge 0}$, and $a\in\mathbb C^\times$. Then:
\begin{enumerate}
\item There exists $b_i\in\{1,2,3\}$ such that, if $m=m_1b_i+m_0$ with $0\le m_0<b_i$ and $T(m\omega_i)$ is the $\lie g[t]$-submodule of $M(b_i\omega_i)^{\otimes m_1}\otimes M(m_0\omega_i)$ generated by the top weight space, then $M(m\omega_i)\cong T(m\omega_i).$
\item $M(m\omega_i)\cong A(m\omega_i)\cong L(\gb\omega_{i,a,m})$.\hfill\qedsymbol
\end{enumerate}
\end{prop}

Our goal is to establish a generalization of the above proposition for minimal affinizations other than Kirillov-Reshetikhin modules. In order to do that, let us introduce the following notation. Given $i\in I, m,r\in \mathbb Z_{\ge 0}$, let $v_{i,m}$ be the image of $1$ in $M(m\omega_i)$ and set
\begin{equation}
R^+(i,m,r) = \{\alpha\in R^+: x_{\alpha,r}^-v_{i,m}=0\}.
\end{equation}
Since $(\lie h\otimes t\mathbb C[t])v_{i,m}=0$, it follows that
\begin{equation}
R^+(i,m,r)\subseteq R^+(i,m,s) \qquad\text{for all }s\ge r.
\end{equation}
The sets $R^+(i,m,r)$ for $\lie g$ not of types $E$ and $F$ were explicitly described in \cite{cha:fer,cm:kr,cm:krg}. We will eventually write them down precisely. For the moment, let us just observe that $R^+(i,m,r)=R^+$ if $r\gg 0$ since $A(m\omega_i)$ is restricted. In fact, if $\lie g$ is of classical type, then $R(i,m,2)=R^+$ for every $i\in I$ and $m\in \mathbb Z_{\ge 0}$. Observe also that $R^+(i,0,0)=R^+$ for all $i\in I$ since $L(0)$ is the trivial representation. Now, given $\lambda\in P^+$ and $r\in\mathbb Z_{\ge 0}$, set
\begin{equation}
R^+(\lambda, r) = \bigcap_{i\in I} R^+(i,\lambda(h_i),r).
\end{equation}
Since $R^+(j,0,s)=R^+$ for all $j\in I$ and $s\in\mathbb Z_{\ge 0}$, it follows that $R^+(m\omega_i,r)=R^+(i,m,r)$ for all $i\in I$ and $m,r\in\mathbb Z_{\ge 0}$ and that \begin{equation}\label{e:Nres}
R^+(\lambda,r)=R^+ \qquad\text{if}\qquad r\gg 0.
\end{equation}

\begin{defn}\label{d:N}
Given $\lambda\in P^+$, let $M(\lambda)$ be the  $\lie g[t]$-module given by the quotient of $U(\lie g[t])$ by the left ideal generated by
\begin{equation}\label{e:Ndef}
\lie n^+[t], \qquad \lie h\otimes t\mathbb C[t], \qquad h_i-\lambda(h_i), \qquad (x_{\alpha_i}^-)^{\lambda(h_i)+1}, \qquad x_{\alpha,r}^-
\end{equation}
for all $i\in I, r\in\mathbb Z_{\ge 0}$, and $\alpha\in R^+(\lambda,r)$.  Let $T(\lambda)$ be the $\lie g[t]$-submodule of $\otm_{i\in I}^{} M(\lambda(h_i)\omega_i)$ generated by the top weight space.
\end{defn}

Definitions \ref{d:KRres} and \ref{d:N} of $M(m\omega_i)$ coincide since $R^+(m\omega_i,r)=R^+(i,m,r)$ for all $i\in I, m,r\in\mathbb Z_{\ge 0}$. The modules $M(\lambda)$ are clearly $\mathbb Z_{\ge 0}$-graded. It follows from Proposition \ref{p:admevrel} that $M(\lambda)$ is a quotient of $A(\lambda)$ and, hence, a restricted $\lie g[t]$-module. Moreover, $T(\lambda)$ is clearly a restricted quotient of $M(\lambda)$ by Proposition \ref{p:cmkr}.

The following is what we expect to be the generalization of Proposition \ref{p:cmkr} when $\lie g$ and $\lambda$ are as in Theorem \ref{t:drimin}.

\begin{con}\label{cj:MN}
Let $\lambda\in P^+$ be such that $\overline{\rm supp}(\lambda)$ does not contain a subdiagram of type $D_4$ and suppose $\gb\lambda\in\cal P_\mathbb A^{++}$ is such that $V_q(\gb\lambda)$ is a minimal affinization of $V_q(\lambda)$. Then, $T(\lambda)\cong M(\lambda)\cong L(\gb\lambda)$.
\end{con}

Proposition \ref{p:cmkr} says the conjecture holds when $\lambda$ is a multiple of a fundamental weight and $\lie g$ is not of type $E$ or $F$. It is quite simple to see that the conjecture also holds when $\lie g$ is of type $A$ for general $\lambda\in P^+$. We now state our main partial results in the direction of proving Conjecture \ref{cj:MN}.

\begin{prop}\label{p:TqM}
Let $\gb\lambda\in\cal P_\mathbb A^{++}$ be such that $V_q(\gb\lambda)$ is a minimal affinization of $V_q(\lambda)$ where $\lambda=\wt(\gb\lambda)$. Then, $T(\lambda)$ is a quotient of $L(\gb\lambda)$.
\end{prop}

\begin{prop}\label{p:MqN}
Let $\lambda\in P^+$ be such that $\overline{\rm supp}(\lambda)$ does not contain a subdiagram of type $D_4$ and suppose $\lie g$ is orthogonal. Then, $L(\gb\lambda)$ is a quotient of $M(\lambda)$.
\end{prop}

\begin{cor}
In the conditions of Proposition \ref{p:MqN}, the first isomorphism of Conjecture \ref{cj:MN} implies the second.\hfill\qedsymbol
\end{cor}

\begin{prop}\label{p:gchar}
Conjecture \ref{cj:MN} holds in the following cases:
\begin{enumerate}
\item $\lie g$ is of type $B$ and ${\rm supp}(\lambda)\subseteq \{1,2,3,n\}$ with $\lambda(h_n)\le 1$ if $n>3$.
\item $\lie g$ is of type $D$ and ${\rm supp}(\lambda)\subseteq (\{1,2,3\}\cap J)\cup\{m\}$ with $m\in\{n-1,n\}$. Here $J=I\backslash \{n-1,n\}$.
\item $\lie g$ is of type $D$ and ${\rm supp}(\lambda)\subseteq \{n-2,n-1,n\}$.
\end{enumerate}
\end{prop}

In the process of proving Proposition \ref{p:gchar} we obtain character formulas for $M(\lambda)$. The proofs of Propositions \ref{p:TqM} and \ref{p:MqN} are given in subsections \ref{ss:tpkr} and \ref{ss:qrel}, respectively. Proposition \ref{p:gchar} is proved in subsections \ref{ss:gcharb} and \ref{ss:gchard}.

\begin{rem}
If $\lie g$ is of classical type, then $R^+(\lambda,2)=R^+$ for every $\lambda\in P^+$ since $R(i,m,2)=R^+$ for every $i\in I,m\in\mathbb Z_{\ge 0}$, as mentioned previously. This implies that the modules $M(\lambda)$ can be regarded as modules for the truncated algebra $\lie g[t]/(\lie g\otimes t^2\mathbb C[t])$ in this case. This was the motivation for the paper \cite{chgr:koszul} where the authors initiated the study of the relations between the finite-dimensional representation theory of $U_q(\tlie g)$ and Koszul algebras. We shall leave the discussion of how our methods are related to those of \cite{chgr:koszul} to a forthcoming publication.
\end{rem}

\section{Tensor products}\label{s:tensor}

\subsection{Tensor products of Kirillov-Reshetikhin modules}\label{ss:tpkr}

The goal of this subsection is to prove Proposition \ref{p:TqM}.
We begin with the following fact which is easily established from \eqref{e:cmLambda}.

\begin{prop}\label{p:dritensor}
Let $\gb\lambda,\gb\mu\in\cal P_q^+$. Then, the $U_q(\tlie g)$-submodule of $V_q(\gb\lambda)\otimes V_q(\gb\mu)$ generated by the top weight space is a quotient of $W_q(\gb\lambda\gb\mu)$.\hfill\qedsymbol
\end{prop}

The following proposition follows from the results of \cite{cha:braid}.

\begin{prop}\label{p:cyclic}
Let $l\in\mathbb Z_{\ge 1}, i_j\in I, m_j\in\mathbb Z_{\ge 1}, a_j\in\mathbb C(q)^\times$ for $j=1,\dots, l$. If $\frac{a_{j}}{a_{k}}\notin q^{\mathbb Z_{> 0}}$ for $j>k$, then $V_q(\gb\omega_{i_1,a_1,m_1})\otimes \cdots\otimes V_q(\gb\omega_{i_l,a_l,m_l})$ is a highest-$\ell$-weight module.\hfill\qedsymbol
\end{prop}

\begin{cor}\label{c:tpkrmin}
Let $\lambda\in P^+, a_i\in\mathbb C(q)^\times, i\in I$, and $\gb\lambda=\prod_{i\in I} \gb\omega_{i,a_i,\lambda(h_i)}$. Then, there exists an ordering $i_1,\dots, i_n$ of $I$ such that $V_q(\gb\lambda)$ is isomorphic to the $U_q(\tlie g)$-submodule of $V_q(\gb\omega_{i_1,a_{i_1},\lambda(h_{i_1})})\otimes\cdots\otimes V_q(\gb\omega_{i_n,a_{i_n},\lambda(h_{i_n})})$ generated by the top weight space.
\end{cor}

\begin{proof}
Let $\gb\omega\in\cal P_q^+$ be such that $\gb\omega^*=\gb\lambda$ and write $\gb\omega_i(u) = \gb\omega_{i,b_i,\lambda^*(h_i)}$ for some $b_i\in\cal P_q^+$. Let also $i'=w_0\cdot i$ for all $i\in I$. It follows from Proposition \ref{p:cyclic} that there exists an ordering $i_1,\dots,i_n$ of $I$ such that $V:=V_q(\gb\omega_{i_n',b_{i_n'},\lambda^*(h_{i_n'})})\otimes\cdots\otimes V_q(\gb\omega_{i_1',b_{i_1'},\lambda^*(h_{i_1'})})$ is highest-$\ell$-weight. Let $W$ be the proper maximal submodule of $V$. Thus, we have a short exact sequence
$$0\to W\to V\to V_q(\gb\omega)\to 0.$$
Then, by Propositions \ref{p:dual} and \ref{p:dualtp}, we also have the following short exact sequence
$$0\to V_q(\gb\lambda)\to V_q(\gb\omega_{i_1,a_{i_1},\lambda(h_{i_1})})\otimes\cdots\otimes V_q(\gb\omega_{i_n,a_{i_n},\lambda(h_{i_n})}) \to W^*\to 0,$$
since $V^*\cong V_q(\gb\omega_{i_1,a_{i_1},\lambda(h_{i_1})})\otimes\cdots\otimes V_q(\gb\omega_{i_n,a_{i_n},\lambda(h_{i_n})})$ and $V_q(\gb\omega)^*\cong V_q(\gb\lambda)$. The corollary now follows immediately.
\end{proof}

We now proceed with the proof of Proposition \ref{p:TqM} as follows.
Given $i\in I$, let $a_i\in\mathbb A^\times$ be such that $\gb\lambda=\prod_{i\in I}\gb\omega_{i,a_i,\lambda(h_i)}$ and let $v_i$ be a highest-$\ell$-weight vector of $V(\gb\omega_{i,a_i,\lambda(h_i)})$. Let also $i_1,\dots,i_n$ be an ordering of $I$ as in Corollary \ref{c:tpkrmin} and $v=v_{i_1}\otimes\cdots\otimes v_{i_n}\in V_q(\gb\omega_{i_1,a_{i_1},\lambda(h_{i_1})})\otimes\cdots\otimes V_q(\gb\omega_{i_n,a_{i_n},\lambda(h_{i_n})})$.
Consider $L_i = U_\mathbb A(\tlie g)v_i, L=U_\mathbb A(\tlie g)v$, and $L'=L_{i_1}\otimes\cdots\otimes L_{i_n}$. Let $a\in\mathbb C^\times$ be such that $\bar{\gb\lambda}=\gb\omega_{\lambda,a}$ and observe that $L(\gb\lambda)\cong \tau_a^*(\overline L)$ and $M(\lambda(h_i)\omega_i)\cong \tau_a^*(\overline{L_i})$, where $\tau_a^*K$ denotes the pullback of a $\lie g[t]$-module $K$ by $\tau_a$. Moreover, it is easy to see that $L\subseteq L'$ and that $\overline{L'}\cong \overline{L_{i_1}}\otimes\cdots\otimes\overline{L_{i_n}}$.

Let $\bar\phi:\overline{L}\to\overline{L'}$ be the map given by Lemma \ref{l:morlim} with $\phi$ being the inclusion $$V_q(\gb\lambda)\to V_q(\gb\omega_{i_1,a_{i_1},\lambda(h_{i_1})})\otimes\cdots\otimes V_q(\gb\omega_{i_n,a_{i_n},\lambda(h_{i_n})}),$$ after identifying $V_q(\gb\lambda)$ with $U_q(\tlie g)v$. It follows that $\tau^*_a(\bar\phi): L(\gb\lambda)\to M(\lambda(h_{i_1})\omega_{i_1})\otimes\cdots\otimes M(\lambda(h_{i_n})\omega_{i_n})$ is a $\lie g[t]$-module map whose image is  $T(\lambda)$.\hfill\qedsymbol

\subsection{A smaller set of relations for $M(\lambda)$}\label{ss:dN}
In this subsection we assume $\lie g$ is orthogonal. Let
$$R^+_1 = \{\alpha\in R^+: \alpha=\sum_{i\in I} n_i\alpha_i \text{ with } n_i\le 1 \text{ for all } i\in I\}.$$
The goal of this subsection is to prove the following proposition.

\begin{prop}\label{p:dN}
For every $\lambda\in P^+$, the module $M(\lambda)$ is isomorphic to the $\lie g[t]$-module $N(\lambda)$ generated by a vector $v$ satisfying
\begin{equation*}
h_iv\ =\ \lambda(h_i)v \qquad\text{and}\qquad \lie n^+[t]v\ =\ \lie h\otimes t\mathbb C[t]v\ =\ (x_{\alpha_i}^-)^{\lambda(h_i)+1}v\ =\ x_{\alpha,1}^-v=0
\end{equation*}
for all $\alpha\in R^+_1$.
\end{prop}

Since every admissible $J\subseteq I$ is of type $A$, it follows that $N(\lambda)$ is a quotient of $A(\lambda)$ and, hence, a finite-dimensional restricted $\lie g[t]$-module. Moreover, it is easy to see that $M(\lambda)$ is a quotient of $N(\lambda)$.
For the converse, set
$$\alpha_{i,j}=\sum_{k=i}^j \alpha_k$$
for all $i,j\in I,i\le j$, if $\lie g$ is of type $B$ and for all $i\le j<n,$ if $\lie g$ is of type $D$.
If $\lie g$ is of type $D$, set also $\alpha_{i,n}=\alpha_{i,n-2}+\alpha_n$ for $i<n-1$ or $i=n$  and $\vartheta_{i}=\alpha_{i,n-1}+\alpha_{n}$ for $i\le n-2$. Furthermore, given $i\le j<n$ ($j<n-2$ if $\lie g$ is of type $D$) set
\begin{alignat*}{3}
&\theta_{i,j} = \alpha_{i,n}+\alpha_{j+1,n} && \text{if } \lie g \text{ is of type } B_n,\\
&\theta_{i,j} = \alpha_{i,n-1}+\alpha_{j+1,n} &\qquad&\text{if } \lie g \text{ is of type } D_n.
\end{alignat*}
Then $R_1^+ = \{\alpha_{i,j}:i,j\in I\}$ ($R^+_1=  \{\alpha_{i,j}:i,j\in I\}\cup \{\vartheta_{i}:i\le n-2\}$ if $\lie g$ is of type $D$) and $R^+=R_1^+\cup \{\theta_{i,j}:i,j\in I\}$.

Denote by $v_{i,m}$ the image of $1$ in $M(m\omega_i), i\in I, m\in\mathbb Z_{\ge 0}$.  Since $R^+(i,m,0)=R^+$ if $m=0$, we shall assume assume $m>0$. Moreover, since we already know that Proposition \ref{p:dN} holds when $\lambda$ is a multiple of a fundamental weight, we assume from now on that $\lambda\in P^+$ is not a multiple of a fundamental weight. From here we split the proof that $N(\lambda)$ is a quotient of $M(\lambda)$ in separate cases according to the type of $\lie g$.

\subsubsection{Type $B$}

It follows from the results of \cite{cha:fer,cm:kr} that
\begin{gather*}
R^+(i,m,0)=\{\alpha_{j,k}: i<j \text{ or } k<i\}\cup \{\theta_{j,k}: i<j\}, \qquad R^+(n,1,1)= R^+(i,m,2)=R^+,\\ \text{for all}\qquad i\in I,\ \ m>0 \quad\text{and}\\
R^+(i,m,1)=R^+(i,m,0)\cup \{\alpha_{j,k}: j\le i\le k\}\cup \{\theta_{j,k}: i\le k\} \quad\text{if}\quad (i,m)\ne(n,1).
\end{gather*}
Set:
\begin{equation}
i_\lambda =
\begin{cases}
\min\{i: \lambda(h_j)=0 \text{ for all } j>i\}, & \text{ if } \lambda(h_n)\ne 1,\\
\min\{i: \lambda(h_j)=0 \text{ for all } i<j<n\}, & \text{ otherwise}.
\end{cases}
\end{equation}
It follows from the above that
\begin{equation}\label{e:Rb}
R^+(\lambda,1)=R^+(i_\lambda,\lambda(h_{i_\lambda}),1)= R^+\backslash\{ \theta_{j,k}: k<i_\lambda\}.
\end{equation}
Proposition \ref{p:dN} follows immediately in the case $\lambda(h_n)>1$.

To complete the proof of Proposition \ref{p:dN}, assume first that $\lambda(h_n)=0$ and notice that $x_{\alpha_{i,j}}^-v=x_{\theta_{i,j}}^-v=0$ if  $i>i_\lambda$. It follows that $x_{\theta_{i,j},r}^-v=[x_{\alpha_{i,n},r}^-,x_{\alpha_{j,n}}^-]v=0$ for all $i,j\in I, j>i_\lambda$ and $r\in\mathbb Z_{>0}$. Also, if $r>1$ and $j\le i_\lambda$, then $x_{\theta_{i,j},r}^-v=[x_{\alpha_{i,n},r-1}^-,x_{\alpha_{j+1,n},1}^-]v=0$. This completes the proof in this case.

If $\lambda(h_n)= 1$, then $x_{\alpha_{i,j}}^-v=0$ if $i_\lambda<i\le j<n$.  Hence, to conclude the proof, it suffices to show that $x_{\theta_{i,n-1},1}^-v=0$ for all $i>i_\lambda$. We prove this inductively on $n-i$. In fact, it follows from the PBW theorem that
$$N(\lambda)[1]\subseteq \sum_{i\le j< n} U(\lie g)x_{\theta_{i,j},1}^-v.$$
In particular, the set of weights of $N(\lambda)[1]$ is contained in $S-Q^+$ where $S=\{\lambda-\theta_{i,j}:i\le j<n\}$.
It is easy to see that $\lambda-\theta_{n-1,n-1}$ is a maximal element of $S$. Hence, if $x_{\theta_{n-1,n-1},1}^-v\ne 0$, it would follow that $V(\lambda-\theta_{n-1,n-1})$ would be an irreducible constituent of $N(\lambda)[1]$. But the condition $\lambda(h_n)=1$ implies $\lambda-\theta_{n-1,n-1}\notin P^+$. Since $N(\lambda)$ is finite-dimensional, it follows that the inductive argument starts. Now suppose we have proved $x_{\theta_{i,n-1},1}^-v=0$ for all $i\ge j$ for some $j\le n-1$ and observe that $\lambda-\theta_{j-1,n-1}$ is a maximal element of $S\backslash\{\lambda-\theta_{i,n-1}:i\ge j\}$. Once more $\lambda-\theta_{j-1,n-1}\notin P^+$ and we conclude the inductive argument as before.

\subsubsection{Type $D$}
In this case we have
\begin{gather*}
R^+(i,m,0)=\{\alpha_{j,k}: i<j \text{ or } k<i\}\cup \{\vartheta_j, \theta_{j,k}: i<j\}\ \text{if}\ i\ne n, n-1,\\
R^+(i,m,0)=\{\alpha_{j,k}: k<n-1 \text{ or } k=i'\} \ \text{if}\ \{i,i'\}=\{n, n-1\},\\
R^+(1,m,1)=R^+(n-1,m,1)=R^+(n,m,1)=R^+(i,m,2)=R^+\ \text{for all}\quad i\in I,\\
R^+(i,m,1)=R^+(i,m,0)\cup \{\alpha_{j,k}: j\le i\le k\}\cup \{\theta_{j,k}: i\le k\} \quad\text{if}\quad i\notin\{1,n-1,n\}.
\end{gather*}
In particular we have
\begin{equation}\label{e:dNdsmall}
R^+(\lambda,1)=R^+ \quad\text{if}\quad \lambda(h_i)=0 \quad\text{for all}\quad i\notin\{1,n-1,n\}
\end{equation}
and, hence, $M(\lambda)$ is an irreducible $\lie g$-module.
Set $i_\lambda = 1$ if $\lambda(h_i)=0$ for all $i\notin\{1,n-1,n\}$ and
\begin{equation}
i_\lambda = \min\{i:  \lambda(h_j)=0 \text{ for all } i<j<n-1\}, \text{ otherwise.}
\end{equation}
It follows  that
\begin{equation}\label{e:Rd}
R^+(\lambda,1)=R^+(i_\lambda,\lambda(h_{i_\lambda}),1)= R^+\backslash\{ \theta_{j,k}: k<i_\lambda\}.
\end{equation}
We are left to show that $x_{\theta_{i,j},1}^-v=0$ if $j>i_\lambda$. But this is clear since $x_{\alpha_{j,n-2}}^-v=0$ if $i_\lambda<j$ and $x_{\theta_{i,j},1}^-=[x_{\alpha_{j,n-2}}^-,x_{\vartheta_{i},1}^-]$.\hfill\qedsymbol

The following corollary is now immediate and proves the first isomorphism of Conjecture \ref{cj:MN} in some very particular cases.

\begin{cor}\label{c:irredN}
Suppose $\lambda\in P^+$ is such that:
\begin{enumerate}
\item $\lambda(h_i)=0$ for all $1<i<n$ and $\lambda(h_n)\le 1$ if $\lie g$ is of type $B$,
\item $\lambda(h_i)=0$ for all $i\notin \{1,n-1,n\}$ if $\lie g$ is of type $D$.
\end{enumerate}
Then, $M(\lambda)$ is irreducible as a $\lie g$-module. In particular, $M(\lambda)\cong T(\lambda)$.\hfill\qedsymbol
\end{cor}

\subsection{The $\ell$-characters}\label{ss:qchar}

Let $\mathbb Z[\cal P_q]$ be the integral group ring over $\cal P_q$. The $\ell$-character of $V\in \wcal C_q$ is defined to be the following element of $\mathbb Z[\cal P_q]$
\begin{equation}
\ch_\ell(V) = \sum_{\gb\mu\in\cal P_q} \dim(V_\gb\mu)\gb\mu.
\end{equation}
The $\ell$-characters are better known as $q$-characters, since this was the name used  when they were first defined in \cite{freer:qchar}. We prefer to call them $\ell$-characters for the following two reasons: first they record information about the dimension of the $\ell$-weight spaces of $V$ (which are not known as $q$-weight spaces), and second, the definition makes sense in the classical context as well. However, due to Theorem \ref{t:tev}, the study of $\ell$-characters in the classical case easily reduces to the study of characters and, therefore, the concept of $\ell$-characters is indeed interesting only in the quantum case.

The proof of the following four results can be found in \cite{cm:chb,freem:qchar}. 

\begin{prop}\label{p:lchasl2}
Let $\lie g=\lie{sl}_2, a\in\mathbb C(q)^\times,$ and $r\in\mathbb Z_{\ge 0}$. Then
$$\ch_\ell(V_q(\gb\omega_{i,a,r}))=\gb\omega_{i,a,r}\sum_{k=0}^r \left(\prod_{j=1}^k \gb\omega_{i,aq^{r-2j},2}\right)^{-1} = \gb\omega_{i,a,r}\sum_{k=0}^r \left(\prod_{j=1}^k \gb\alpha_{i,aq^{r-2j+1}}\right)^{-1}.$$
\hfill\qedsymbol
\end{prop}

\begin{thm}\label{t:cone}
Let $V$ be a quotient of $W_q(\gb\lambda)$ for some $\gb\lambda\in\cal P_q^+$. If $V_\gb\mu\ne 0$, then $\gb\mu\le\gb\lambda$.\hfill\qedsymbol
\end{thm}

\begin{prop}\label{p:cone}
Let $V\in\wcal C_q, v\in V_\gb\mu\backslash\{0\}$ for some $\gb\mu\in\cal P_q$, and suppose $i\in I$ is such that $x_{i,r}^+v=0$ for all $r\in\mathbb Z$. Then, $\gb\mu_i(u)\in\mathbb C(q)[u]$ and, if $\gb\mu_i(u) = \prod_{k=1}^m f_{i,a_k,r_k}(u)$ as in \eqref{e:qsf},
$$x_i^-v\in \sum_{k=1}^{m}\sum_{j=1}^{r_k} V_{\gb\mu(\gb\alpha_{i,a_kq_i^{r_k+1-2j}})^{-1}}.$$
Moreover, $\dim(V_{\gb\mu(\gb\alpha_{i,a_kq^{r_k-1}})^{-1}})\ge \#\{1\le l\le k:a_l=a_k\}.$
\hfill\qedsymbol
\end{prop}

Given $V\in\wcal C_q$, let
$$\wtl(V) = \{\gb\mu\in\cal P_q: V_\gb\mu\ne 0\}.$$
A highest-$\ell$-weight module $V$ of highest $\ell$-weight $\gb\lambda\in\cal P_q^+$ is said to be special if
$$\wtl(V)\cap\cal P_q^+ = \{\gb\lambda\}.$$

\begin{thm}\label{t:fmas}
If $\gb\lambda\in\cal P_q^+$ is such that $V_q(\gb\lambda)$ is special, then the output of the Frenkel-Mukhin algorithm with input $\gb\lambda$ is $\chl(V_q(\gb\lambda))$.\hfill\qedsymbol
\end{thm}

The following theorem was proved in \cite{her:min}.

\begin{thm}
If $\lie g$ is of type $A,B$, or $G$, then all minimal affinizations are special.\hfill\qedsymbol
\end{thm}

Let $\lie g$ be of type $A, B$, or $G$, and $\lambda\in P^+$. It follows from the above that if $V_q(\gb\lambda)$ is a minimal affinization of $V_q(\lambda)$, then  $\chl(V_q(\gb\lambda))$ is given by the Frenkel-Mukhin algorithm. We will actually need only the following corollary of the algorithm. Let $V\in\wcal C_q, v\in V_\gb\mu\backslash\{0\}$ for some $\gb\mu\in\cal P_q$, and suppose $i\in I$ is such that $x_{i,r}^+v=0$ for all $r\in\mathbb Z$. Using Proposition \ref{p:cone}, we can write $\gb\mu_i(u) = \prod_{k=1}^m f_{i,a_k,r_k}(u)$  as in \eqref{e:qsf}. Then, the algorithm implies that
\begin{equation}\label{e:fma1}
\gb\mu\gb\alpha_{i,b}^{-1}\in \wtl(V_q(\gb\lambda)) \qquad\text{iff}\qquad b = a_kq_i^{r_k-1} \qquad\text{for some}\qquad k=1,\dots,m.
\end{equation}
The next proposition will be crucial for the proof of Proposition \ref{p:MqN}.

\begin{prop}\label{p:lchA}
Suppose $\lie g$ is of type $A$, $\lambda\in P^+$, $\gb\lambda=\prod_{i\in I}\gb\omega_{i,a_i,\lambda(h_i)}$, $\gb\mu\in\wtl(V_q(\gb\lambda))$, and  $\gb\lambda\gb\mu^{-1} = \gb\alpha_{j,b_j}\gb\alpha_{j+1,b_{j+1}}\cdots\gb\alpha_{k,b_k}$ for some $j\le k$ and some $a_i,b_l\in\mathbb C(q)^\times, i\in I, l=j,\dots,k$.
\begin{enumerate}
\item If $\frac{a_{i+1}}{a_i}=q^{\lambda(h_i)+\lambda(h_{i+1})+1}$ for all $i<n$, then $b_k=a_kq^{\lambda(h_k)-1}$.
\item If $\frac{a_{i+1}}{a_i}=q^{-(\lambda(h_i)+\lambda(h_{i+1})+1)}$ for all $i<n$, then $b_j=a_jq^{\lambda(h_j)-1}$.
\end{enumerate}
\end{prop}

\begin{proof}
Straightforward using induction on $k-j$ together with \eqref{e:fma1}.
\end{proof}

\subsection{Quantized relations}\label{ss:qrel}
We now prove Proposition \ref{p:MqN}. In particular, we assume that $\lie g$ is orthogonal. To make the notation more uniform, we assume for the rest of the proof that $\lie g$ is of type $B_n$ or $D_{n+1}, n\ge 2$. Before we begin, let us remark the following corollary of Proposition \ref{p:MqN} and Corollary \ref{c:irredN}.

\begin{cor}
Conjecture \ref{cj:MN} holds if $\lambda$ satisfies the conditions of Proposition \ref{p:MqN} as well as of Corollary \ref{c:irredN}. In particular, if $V_q(\gb\lambda)$ is a minimal affinization of $V_q(\lambda)$, then  $V_q(\gb\lambda)\cong V_q(\lambda)$ as a $U_q(\lie g)$-module.\hfill\qedsymbol
\end{cor}

If $\lambda$ is supported on an admissible subdiagram, Proposition \ref{p:MqN} easily follows from Propositions \ref{p:admevrel} and \ref{p:dN}. In particular, we can henceforth assume that the support of $\lambda$ contains a spin node  and that there exists $i<n$ such that $\lambda(h_i)\ne 0$.  If $\lie g$ is of type $D$, we will prove Proposition \ref{p:MqN} in the case $\lambda(h_{n+1})=0$ (the other cases are proved similarly). Set
\begin{equation}
i_\lambda = \min\{i:  \lambda(h_j)=0 \text{ for all } i<j<n\}, \text{ otherwise.}
\end{equation}
Observe that the above definition of $i_\lambda$ does not coincide with the one given in subsection \ref{ss:dN} for $\lie g$  of type $B$ and $\lambda(h_n)>1$.

From now on we assume that $\gb\lambda\in\cal P_\mathbb A^{++}$ is such that $V=V_q(\gb\lambda)$ is a minimal affinization of $V_q(\lambda)$ and $a\in\mathbb C$ is such that $\bar{\gb\lambda}=\gb\omega_{\lambda,a}$. We also fix a highest-$\ell$-weight $v$ vector of $V$ and $a_i\in\mathbb A^\times, i\in I$, such that $\gb\lambda=\prod_{i\in I} \gb\omega_{i,a_i,\lambda(h_i)}$.

Let $v'$ be the image of $v$ in $L(\gb\lambda)$.
It again follows from Proposition \ref{p:admevrel} that $x_{\alpha_{i,j},1}^-v'=0$ if $i=j$ or if $j<n$. If $\lie g$ is of type $B$, this implies $x_{\alpha_{i,n},1}^-v'=0$ if $i>i_\lambda$. If $\lie g$ is of type $D$ and $i_\lambda<i<n$, it follows that $x_{\alpha_{i,n},1}^-v'=x_{\alpha_{i,n+1},1}^-v'=x_{\vartheta_{i},1}^-v'=0$. We claim that it remains to show that $x_{\alpha_{i,n},1}^-v'=0$ for $i\le i_\lambda$. In fact, if $\lie g$ is of type $B$ this is clear from  Proposition \ref{p:dN}. If $\lie g$ is of type $D$, it follows that $x_{\vartheta_i,1}^-v' = [x_{\alpha_{n+1}}^-,x_{\alpha_{i,n},1}^-]v'=0$ provided $x_{\alpha_{i,n},1}^-v'=0$.

Denote by $\bar v$ the image of $v$ in $\overline V$. It suffices to show that $x_{\alpha_{i,n},1}^-\bar v = ax_{\alpha_{i,n}}^-\bar v$ for all $i\le i_\lambda$. Consider the elements
\begin{equation}
X_{\alpha_{i,j},r}^- = [x_{j,r}^-,[x_{j-1}^-,\dots[x_{i+1}^-,x_i^-]\dots]] \qquad\text{and}\qquad k_{i,j} = \prod_{l=i}^j k_l
\end{equation}
for $i\le j\le n$ and $r\in\mathbb Z_{\ge 0}$. Set also $X_{\alpha_{i,j}}^-:=X_{\alpha_{i,j},0}^-$ and notice that $X_{\alpha_{i,j},r}^- = [x_{j,r}^-,X_{\alpha_{i,j-1}}^-]$ if $i<j$.
Clearly $X_{\alpha_{i,j},r}^-\in U_\mathbb A(\tlie n^-)$ and $\overline{X_{\alpha_{i,j},r}^-}=x_{\alpha_{i,j},r}^-$. We will need the following two lemmas
\begin{lem}\label{l:MqNo}
Suppose $V\in\cal C_q, \mu\in P,w\in V_\mu\backslash\{0\}$, and $i\in I$ are such that $\mu(h_l)=0$ and $x_{l}^+w=0$ for $l> i$. Then $X_{\alpha_{i,j}}^-w = x_{j}^-x_{j-1}^-\cdots x_i^-w$.
\end{lem}

\begin{proof}
It is a straightforward computation using the commutation relations $[x_{l,s}^-,x_{l',s'}^-]=0$ for $l,l'$ such that $c_{l,l'}=0$ and $x_l^-w=0$ if $l>i$.
\end{proof}

The next lemma follows from the relation $[h_{i,1},x_i^-]=-[2]_{q_i}x_{i,1}^-$ together with Proposition \ref{p:lchasl2}.

\begin{lem}\label{l:x1}
Suppose  $w$ is a highest-$\ell$-weight vector of $V_q(\gb\omega_{i,a_i,m})$ for some $i\in I,$ and some $m\in\mathbb Z_{\ge 0}$, then $x_{i,1}^-w = a_iq_i^{m}x_i^-w.$\hfill\qedsymbol
\end{lem}

Let $\gb\lambda'$ be such that $\gb\lambda=\gb\lambda'\gb\omega_{n,a_n,\lambda(h_n)}$. Let also $v_1,v_2$ be highest-$\ell$-weight vectors of $V_q(\gb\lambda')$ and $V_q(\gb\omega_{n,a_n,\lambda(h_n)})$.
By Proposition \ref{p:cyclic} and Corollary \ref{c:tpkrmin}, either $V_q(\gb\lambda)\cong U_q(\tlie g)(v_1\otimes v_2)\subseteq V_q(\gb\lambda')\otimes V_q(\gb\omega_{n,a_n,\lambda(h_n)})$ or $V_q(\gb\lambda)\cong U_q(\tlie g)(v_2\otimes v_1)\subseteq  V_q(\gb\omega_{n,a_n,\lambda(h_n)})\otimes V_q(\gb\lambda')$. We assume we are in the former case (the latter is proved similarly using part (b) of Proposition \ref{p:lchA} instead of part (a)). In particular, by Theorem \ref{t:drimin},  we must have
\begin{equation}\label{e:minpos}
\frac{a_{i+1}}{a_i} = q^{d_i\lambda(h_i)+d_{i+1}\lambda(h_{i+1})+r^\vee} \qquad\text{for all}\qquad i<n.
\end{equation}

By Lemmas \ref{l:comcom} and \ref{l:comult}, modulo elements of the form $x(v_1\otimes v_2)$ with $x\in U_\mathbb A(\tlie g)\otimes U_\mathbb A(\tlie g)$ such that $\bar x=0$, we have
\begin{align*}
X_{\alpha_{i,n}}^-(v_1\otimes v_2) & = x_{n}^-X_{\alpha_{i,n-1}}^-(v_1\otimes v_2)- X_{\alpha_{i,n-1}}^-x_{n}^-(v_1\otimes v_2)\\
& = x_n^-((X_{\alpha_{i,n}}^-v_1)\otimes v_2)-X_{\alpha_{i,n-1}}^-(v_1\otimes (x_n^-v_2))\\
& = (x_{n}^-X_{\alpha_{i,n-1}}^-v_1)\otimes (q^{-\lambda(h_n)}  v_2) + (X_{\alpha_{i,n-1}}^-v_1)\otimes (x_{n}^- v_2)\\
& - (X_{\alpha_{i,n-1}}^-v_1)\otimes (k_{i,n-1}^{-1}x_{n}^-v_2)- v_1\otimes(X_{\alpha_{i,n-1}}^-x_n^-v_2)\\
& = q^{-\lambda(h_n)}(x_n^-X_{\alpha_{i,n-1}}^-v_1)\otimes v_2 + (1-q^{-2})(X_{\alpha_{i,n-1}}^-v_1)\otimes(x_n^-v_2) - v_1\otimes(X_{\alpha_{i,n-1}}^-x_n^-v_2).
\end{align*}
On the other hand,
\begin{align*}
X_{\alpha_{i,n},1}^-(v_1\otimes v_2) & = x_{n,1}^-X_{\alpha_{i,n-1}}^-(v_1\otimes v_2)- X_{\alpha_{i,n-1}}^-x_{n,1}^-(v_1\otimes v_2)\\
& = x_{n,1}^-((X_{\alpha_{i,n}}^-v_1)\otimes v_2)-X_{\alpha_{i,n-1}}^-(v_1\otimes (x_{n,1}^-v_2))\\
& = (x_{n,1}^-X_{\alpha_{i,n-1}}^-v_1)\otimes (q^{\lambda(h_n)}  v_2) + (X_{\alpha_{i,n-1}}^-v_1)\otimes (x_{n,1}^- v_2)\\
& - (X_{\alpha_{i,n-1}}^-v_1)\otimes (k_{i,n-1}^{-1}x_{n,1}^-v_2)- v_1\otimes(X_{\alpha_{i,n-1}}^-x_{n,1}^-v_2)\\
& = q^{\lambda(h_n)}(x_{n,1}^-X_{\alpha_{i,n-1}}^-v_1)\otimes v_2 + (1-q^{-2})(X_{\alpha_{i,n-1}}^-v_1)\otimes(x_{n,1}^-v_2) - v_1\otimes(X_{\alpha_{i,n-1}}^-x_{n,1}^-v_2).
\end{align*}
Using Lemma \ref{l:x1} we get
\begin{align*}
X_{\alpha_{i,n},1}^-(v_1\otimes v_2) & = q^{\lambda(h_n)}(x_{n,1}^-X_{\alpha_{i,n-1}}^-v_1)\otimes v_2\\
& +a_nq^{\lambda(h_n)}\left((1-q^{-2})(X_{\alpha_{i,n-1}}^-v_1)\otimes(x_{n}^-v_2) - v_1\otimes(X_{\alpha_{i,n-1}}^-x_{n}^-v_2)\right)\\
&= a_nq^{\lambda(h_n)}X_{\alpha_{i,n}}^-(v_1\otimes v_2) + q^{\lambda(h_n)}(x_{n,1}^-X_{\alpha_{i,n-1}}^-v_1)\otimes v_2 - a_n (x_n^-X_{\alpha_{i,n-1}}^-v_1)\otimes v_2.
\end{align*}
Hence, it suffices to show that
\begin{equation}\label{e:MqNend}
q^{\lambda(h_n)}(x_{n,1}^-X_{\alpha_{i,n-1}}^-v_1)\otimes v_2 = a_n (x_n^-X_{\alpha_{i,n-1}}^-v_1)\otimes v_2.
\end{equation}
If $i>i_\lambda$, both sides of the above equality vanish. If $i\le i_\lambda$ we proceed as follows.
Notice that $x_{n,r}^+X_{\alpha_{i,n-1}}^-v_1=0$ for all $r\in\mathbb Z$ and let $W$ be the $U_q(\tlie g_n)$-submodule of $V_q(\gb\lambda')$ generated by $X_{\alpha_{i,n-1}}^-v_1$. Then, by Proposition \ref{p:lchA}(a), the highest-$\ell$-weight of $W$ is $\gb\omega_{n,a_{n-1}q^{r^\vee\lambda(h_{n-1})},r^\vee}$. Moreover, by \eqref{e:fma1}, $W$ is a minimal affinization. Hence, by Lemma \ref{l:x1},
$$x_{n,1}^-X_{\alpha_{i,n-1}}^-v_1 = a_{n-1}q^{r^\vee(\lambda(h_{n-1})+1)}x_{n}^-X_{\alpha_{i,n-1}}^-v_1.$$
This and \eqref{e:minpos} imply \eqref{e:MqNend}.\hfill\qedsymbol

\section{Graded Characters of Restricted Limits of Minimal Affinizations}\label{s:char}

\subsection{Preliminaries}

Although Theorem \ref{t:drimin} tells which objects of $\wcal C_q$ correspond to minimal affinizations, it does not say anything about their $U_q(\lie g)$-structure, unless $\lie g$ is of type $A$. In some few cases this is known (see \cite{cha:minr2,cm:kr,cm:krg}).  Naturally, in principle, the $U_q(\lie g)$-structure can be read off the $\ell$-character. In practice, this is not so easy to do, even in the situations that the Frenkel-Mukhin algorithm does produced the $\ell$-character. We will now apply the techniques of \cite{cm:kr,cm:krg} to prove Proposition \ref{p:gchar} and, hence, Conjecture \ref{cj:MN} in those cases. As a byproduct of the proof, we obtain closed formulas for the character of the minimal affinizations if $\lambda$ is as in Proposition \ref{p:gchar}. We shall also prove an analogue of Conjecture \ref{cj:MN} in the case of multiple minimal affinizations for $\lie g$ of type $D_4$.

We shall make use of the following lemma (see \cite[\S1.5]{cm:krg}).

\begin{lem}\label{l:heis}
Consider the three-dimensional Heisenberg algebra $\mathfrak H$ spanned by elements $x,y,z$ where $z$ is central and $[x,y]=z$. Suppose that $V$ is a  representation of $\frak H$ and let $0\ne v\in V$ be such that $x^rv=0$. Then  for all $k,s\in\mathbb Z_{\ge 0}$  the element $y^kz^sv$ is in the span of elements of the form $x^ay^bz^cv$ with  $0\le c<r, a+c=s, b+c=k+s$.\hfill\qedsymbol
\end{lem}

Introduce the following notation. Let $V$ be a finite-dimensional $\lie g$-module and $\lambda\in P^+$. Denote by $m_\lambda(V)$ the multiplicity of the irreducible module $V(\lambda)$ as an irreducible constituent of $V$.
Set $d_i' = d_i/r^\vee$. Hence, if $\lie g$ is simply laced, $d_i'=d_i=1$. If $\lie g$ is of type $B, d_i'=1$ if $i<n$ and $d_n'=\frac{1}{2}$. The symbol $[m]$ means the largest integer not greater than $m$.

\subsection{Type $B$}\label{ss:gcharb}

Given  $i\in I$ and $m\in\mathbb Z_{\ge 0}$, let $v_{i,m}$ be the image of $1$ in $M(m\omega_i)$. The following was proved in \cite{cha:fer,cm:kr}.

\begin{lem}\label{l:krb} \hfill
\begin{enumerate}
\item $M(m\omega_1)\cong V(m\omega_1)$.
\item $M(m\omega_2)[l]=0$ if $l>[d_2'm]$ and $M(m\omega_2)[l] = U(\lie n^-)(x_{\theta_{1,1},1}^-)^lv_{2,m}\cong V(m\omega_2-l\theta_{1,1})$ if $0\le l\le [d_2'm]$.
\item $M(m\omega_3)[l]=0$ if $l>[d_3'm]$ and $M(m\omega_3)[l] = U(\lie n^-)(x_{\theta_{2,2},1}^-)^lv_{3,m}\cong V(m\omega_3-l\theta_{2,2})$ if $0\le l\le [d_3'm]$. Moreover, $(x_{\theta_{1,2},1}^-)^{r_{1,2}}(x_{\theta_{2,2},1}^-)^{r_{2,2}}v_{3,m}$ is a multiple of $(x_{\alpha_1}^-)^{r_{1,2}}(x_{\theta_{2,2},1}^-)^{r_{1,2}+r_{2,2}}v_{3,m}$.\hfill\qedsymbol
\end{enumerate}
\end{lem}

The ``moreover'' part of the above lemma can also be proved using Lemma \ref{l:heis}.

\begin{prop}
Suppose $\lambda\in P^+$ is such that $\lambda(h_i)=0$ for $i>2$. Then,   $M(\lambda)[k]\cong  V(\lambda-k\theta_{1,1})$ if $0\le k\le [d_2'\lambda(h_2)]$ and $M(\lambda)[k]=0$ otherwise. Moreover, $M(\lambda)\cong T(\lambda)$.
\end{prop}

\begin{proof}
Let $v$ be the image of $1$ in $M(\lambda)$.  Equation \eqref{e:Rb} implies $R^+(\lambda,1)\supseteq R^+\backslash\{\theta_{1,1}\}$. Together with the PBW Theorem, this implies that
$$M(\lambda) = \opl_{k\ge 0}^{} U(\lie n^-)(x_{\theta_{1,1},1}^-)^kv.$$
It follows from Lemma \ref{l:hwvecs} that $m_\mu(M(\lambda))\le 1$ for every $\mu$ and that $m_\mu(M(\lambda))$ may be nonzero only when $\mu=\lambda-k\theta_{1,1}$ for some $k\in\mathbb Z_{\ge 0}$. Since $\theta_{1,1} = (d_2')^{-1}\omega_2$, $\mu-k\theta_{1,1}\in P^+$ iff $k\le [d_2'\lambda(h_2)]$.

Now let $v_i$ be a nonzero element in the top weight space of  $M(\lambda(h_i)\omega_i), i=1,2$. Then, by \cite{cm:kr}, $(x_{\theta_{1,1},1}^-)^kv_2$ is the highest-weight vector of the irreducible $\lie g$-submodule of $M(\lambda(h_2)\omega_2)[k], k=0,\dots,[d_2'\lambda(h_2)]$, while $M(\lambda(h_1)\omega_1)$ is an irreducible  $\lie g$-module itself. Therefore, $(x_{\theta_{1,1},1}^-)^k(v_1\otimes v_2)=v_1\otimes (x_{\theta_{1,1},1}^-)^kv_2\ne 0$ proving that $T(\lambda)[k]\ne 0$. Hence, $T(\lambda)[k]\cong M(\lambda)[k]$ and we are done.
\end{proof}

\begin{rem}
In particular, the above proposition reproves one of the the main results of \cite{cha:minr2} using a different method.
\end{rem}

Now assume $n\ge 3$ and suppose $\lambda\in P^+$ is such that $\lambda(h_i)=0$ for $i>3$. In this case equation \eqref{e:Rb} implies
\begin{equation}\label{e:Rb3}
R^+(\lambda,1)\supseteq R^+\backslash\{\theta_{2,2},\theta_{1,2},\theta_{1,1}\}.
\end{equation}
Observe that $\theta_{1,1}=\omega_2, \theta_{1,2}=\omega_1-\omega_2+(d_3')^{-1}\omega_3$ and $\theta_{2,2}=(d_3')^{-1}\omega_3-\omega_1$. In particular, $\{\theta_{2,2},\theta_{1,2},\theta_{1,1}\}$ is a linearly independent subset of $\lie h^*$. Let $\gbr e_j, j\in\mathbb Z_{\ge 0}$, be the standard basis of $\mathbb Z_{\ge 0}^3$, set
\begin{equation}\label{e:b3pset}
\cal A_3(\lambda) = \left\{\gbr r=(r_1,r_2,r_3)\in\mathbb Z_{\ge 0}^3: r_3\le\lambda(h_2), r_2\le\lambda(h_1),r_1+r_2\le [d_3'\lambda(h_3)]\right\}
\end{equation}
and, given $\gb r\in \mathbb Z_{\ge 0}^3$, define
\begin{equation}
\gbr y_\gbr r = (x_{\theta_{2,2},1}^-)^{r_1}(x_{\theta_{1,2},1}^-)^{r_2}(x_{\theta_{1,1},1}^-)^{r_3}.
\end{equation}
Notice that the elements $x_{\theta_{2,2},1}^-,x_{\theta_{1,2},1}^-,x_{\theta_{1,1},1}^-$ commute among themselves.

\begin{lem}\label{l:ubb3}
Let $v$ be the image of $1$ in $M(\lambda)$. For every $\gbr s\in\mathbb Z_{\ge 0}^3$,
$$\gbr y_\gbr sv\in\sum_{\gbr r} U(\lie n^-)\gbr y_\gbr rv$$
where the sum is over the elements $\gbr r\in\mathbb Z_{\ge 0}^3$ such that $r_3\le\lambda(h_2)$ and $r_2\le\lambda(h_1)$.
\end{lem}

\begin{proof}
By Lemma \ref{l:heis} with $x=x_{\alpha_2}^-, y=x_{\theta_{1,2},1}^-, z=x_{\theta_{1,1},1}^-$ we have that $\gbr y_\gbr sv$ is in the span of elements of the form $(x_{\alpha_2}^-)^a\gbr y_{\gbr s'}v$ with $a>0$ and $\gbr s'$ such that $s_3'\le \lambda(h_2)$. Using Lemma \ref{l:heis} once more, this time with $x=x_{\alpha_1}^-, y=x_{\theta_{2,2},1}^-, z=x_{\theta_{1,2},1}^-$, it follows that an element $\gbr y_{\gbr s'}v$ with $\gbr s'$ as above belongs to the span of elements of the form $(x_{\alpha_1}^-)^a\gbr y_{\gbr r}v$ with $a>0$ and $\gbr r$ as claimed.
\end{proof}

Given $\gbr r\in\mathbb Z_{\ge 0}^3$, define
\begin{equation}\label{e:wt(r)b}
\wt(\gbr r) = r_1\theta_{2,2}+r_2\theta_{1,2}+r_3\theta_{1,1}\qquad\text{and}\qquad \gr(\gbr r) = r_1+r_2+r_3.
\end{equation}
Since $\{\theta_{2,2},\theta_{1,2},\theta_{1,1}\}$ is linearly independent, it follows that $\wt$ is an injective function.

\begin{prop}\label{p:Nb3}
For every $\lambda\in P^+$ as above we have $M(\lambda)\cong T(\lambda)$ and
$$M(\lambda)[l]\cong \opl_{\gbr r\in\cal A_3(\lambda): \gr(\gbr r)=l}^{} V(\lambda-\wt(\gbr r)).$$
\end{prop}

\begin{proof}
Let $v$ be the image of $1$ in $M(\lambda)$.  Equation \eqref{e:Rb3}, together with the PBW Theorem, implies that
$$M(\lambda) = \sum_{\gbr r\in\mathbb Z_{\ge 0}^3}^{} U(\lie n^-)\gbr y_\gbr rv.$$
Lemma \ref{l:ubb3} implies that the above sum can be restricted to $\gbr r$ such that $r_3\le\lambda(h_2)$ and $r_2\le\lambda(h_1)$. This, together with Lemma \ref{l:hwvecs}, implies that $m_\mu(M(\lambda))\le 1$  and equality may occur only if $\mu=\lambda- \wt(\gb r)$ for some $\gbr r$ as above. Moreover, $\wt(\gbr r)\in P^+$  only if $r_1+r_2\le [d_3'\lambda(h_3)]$ and, hence, $\gbr r$ must be in $\cal A_3(\lambda)$. It follows that $M(\lambda)[l]$ is a quotient of $\opl_{\gbr r\in\cal A_3(\lambda): \gr(\gbr r)=l}^{} V(\lambda-\wt(\gbr r))$. In order to complete the proof, it suffices to show that $T(\lambda)[l]$ contains a submodule isomorphic to $V(\lambda-\wt(\gbr r))$ for every $\gbr r\in\cal A_3(\lambda)$ such that $\gr(\gbr r)=l$.

Thus, let $v_i=v_{i,\lambda(h_i)}, i=1,2,3$, and let $v_i^j$ be the image of $v_i$ in $M(\lambda(h_i)\omega_i)(j)$ for $j\ge 0$.
Then, if $\gbr r\in\cal A_3(\lambda)$, Lemma \ref{l:krb} implies
\begin{align}\notag
\gbr y_\gbr r(v_1\otimes v_2\otimes v_3^{r_1+r_2})= &\ v_1\otimes (x_{\theta_{1,1},1}^-)^{r_3} v_2\otimes (x_{\theta_{2,2},1}^-)^{r_1}(x_{\theta_{1,2},1}^-)^{r_2}v_3^{r_1+r_2}\\ \label{e:Nb3}\\\notag
 = &\ v_1\otimes (x_{\theta_{1,1},1}^-)^{r_3} v_2\otimes (x_{\alpha_1}^-)^{r_2}(x_{\theta_{2,2},1}^-)^{r_1+r_2}v_3^{r_1+r_2}\ne 0.
\end{align}
Given $r\le\lambda(h_2), s\le [d_3'\lambda(h_3)]$, let $T_{r,s}(\lambda)$ be the $\lie g[t]$-submodule of $M(\lambda(h_1)\omega_1)\otimes M(\lambda(h_2)\omega_2)(r)\otimes M(\lambda(h_3)\omega_3)(s)$ generated by $v_{r,s}:=v_1\otimes v_2^r\otimes v_3^s$. Clearly $T_{r,s}(\lambda)$ is a quotient of $T(\lambda)(r+s)$.
Set $\gbr r_0=(s,0,r)$ and $\gbr r_j = \gbr r_{j-1}+(\gbr e_2-\gbr e_1)$ for $1\le j\le s':=\min(\lambda(h_1),s)$.
Notice that
$$T_{r,s}(\lambda)[r+s]=\sum_{j=0}^{s'} U(\lie n^-)\gbr y_{\gbr r_j}v_{r,s} \qquad\text{and}\qquad (\lambda -\wt(\gbr r_j))(h_1)=\lambda(h_1)+s.$$
In particular, $\lambda-\wt(\gbr r_0)$ is the unique maximal weight of $T_{r,s}(\lambda)[r+s]$.
We prove inductively on $k=0,1,\dots,s'$ that
$$\sum_{j=0}^k U(\lie n^-)\gbr y_{\gbr r_j}v_{r,s} \cong \opl_{j=0}^k V(\lambda-\wt(\gbr r_j))$$
as $\lie g$-module. Since every $\gbr r\in\cal A_3(\lambda)$ is of the form $\gbr r_j$ for some $r,s,j$ as above this completes the proof.

It is clear from Lemma \ref{l:krb} and \eqref{e:Nb3} that $\lie n^+\gbr y_{\gbr r_0}v_{r,s}=0$ and, hence, generates a $\lie g$-submodule isomorphic to $V(\lambda-\wt(\gbr r_0))$. In particular, we can assume $s'>0$. Notice that the weight space of $V(\lambda-\wt(\gbr r_j))$ of weight $\lambda-\wt(\gbr r_j)-(k-j)\alpha_1$ is one-dimensional for $0\le j\le k$. Using the induction hypothesis on $k$, we know that the dimension of the weight space of $\sum_{j=0}^k U(\lie n^-)\gbr y_{\gbr r_j}v_{r,s}$ of weight $\lambda-\wt(\gbr r_0)-(k+1)\alpha_1$ is $k+1$. Since the elements $(x_{\alpha_1}^-)^{j}\gbr y_{\gbr r_{k+1-j}}v_{r,s}, 0\le j\le k+1$ are clearly linearly independent, it follows that $V(\lambda-\wt(\gbr r_{k+1}))$ is a submodule of $\sum_{j=0}^{k+1} U(\lie n^-)\gbr y_{\gbr r_j}v_{r,s}$.
\end{proof}

\begin{rem}
Suppose $n>3$ and that $\lambda\in P^+$ is such that ${\rm supp}(\lambda)\subseteq\{1,2,3,n\}$ with $\lambda(h_n)=1$. Since $R^+(n,1,1)=R^+$, it follows that all of the above can be carried out and Proposition \ref{p:Nb3} remains valid (notice $d_3'=1$ in this case).
\end{rem}

\subsection{Type $D$}\label{ss:gchard}  Define the set $\cal A_3(\lambda)$ exactly as in \eqref{e:b3pset} and the maps $\wt$ and ${\rm gr}$ as in \eqref{e:wt(r)b}.

\begin{prop}\label{p:nd3}
If $\lambda\in P^+$ is such that $\lambda(h_i)=0$ if $3<i<n-1$  then, $M(\lambda)\cong T(\lambda)$.  Moreover:
\begin{enumerate}
\item If $n=4$, then $M(\lambda)[l]\cong  V(\lambda-l\theta_{1,1})$ if $0\le l\le \lambda(h_2)$ and $M(\lambda)[l]=0$ otherwise.
\item If $n>4$, then
$$M(\lambda)[l]\cong \opl_{\gbr r\in\cal A_3(\lambda): \gr(\gbr r)=l}^{} V(\lambda-\wt(\gbr r)).$$
\end{enumerate}
\end{prop}

\begin{proof}
The proof is essentially the same as that of Proposition \ref{p:Nb3} using that $R^+(m\omega_i,1)=R^+$ if $i$ labels a spin node.
\end{proof}

In particular, the above proposition gives the description of the graded characters of $M(\lambda)$ in types $D_4$ and $D_5$ for any $\lambda\in P^+$. If ${\rm supp}(\lambda)$ contains at most one of the spin nodes, it follows that the above is also the character of the minimal affinizations of $V_q(\lambda)$. Otherwise, it is just a lower bound.

Now, let $m\in\mathbb Z_{\ge 0}$ and set $\cal A(m)=\{\gbr r=(r_1,r_2,\dots,r_{[(n-2)/2]})\in\mathbb Z_{\ge 0}^{[(n-2)/2]}: m\ge r_1\ge\cdots\ge r_{[(n-2)/2]}\}$. Define
\begin{equation}
\wt(\gbr r) = \sum_{j=1}^{[(n-2)/2]} r_j\theta_{n-2j,n-2j}, \qquad\text{and}\qquad \gr(\gbr r) = \sum_{j=1}^{[(n-2)/2]} r_j.
\end{equation}
It was proved in \cite{cha:fer,cm:kr} that
\begin{equation}
M(m\omega_{n-2})[l] = \sum_{\gbr r\in\cal A(m):\gr(\gbr r)=l} \gbr y_\gbr rv_{n-2,m}\cong \opl_{\gbr r\in\cal A(m):\gr(\gbr r)=l}^{} V(m\omega_{n-2}-\wt(\gbr r))
\end{equation}
where
\begin{equation}
\gbr y_\gbr r =\prod_{j=1}^{[(n-2)/2]} (x_{\theta_{n-2j,n-2j},1}^-)^{r_j}.
\end{equation}
Proceeding similarly to the proof of Proposition \ref{p:Nb3} one also proves the following.
\begin{prop}\label{p:onlyspin}
Let $\lambda\in P^+$ be such that $\lambda(h_i)=0$ if $i<n-2$.  Then, $M(\lambda)\cong T(\lambda)$ and
\begin{equation*}
M(\lambda)[l] \cong \opl_{\gbr r\in\cal A(\lambda(h_{n-2})): \gr(\gbr r)=l}^{} V(\lambda-\wt(\gbr r)).
\end{equation*}
\hfill\qedsymbol
\end{prop}

\subsection{Multiple minimal affinizations: the regular case}\label{ss:multmar}

Let $\lie g$ be of types $D$ or $E$ and $i_0\in I$ be the unique node triply connected. Let also $J_1, J_2, J_3\subseteq I$ be an enumeration of the three maximal subdiagrams of type $A$ of the Dynkin diagram of $\lie g$ (they are not admissible). Let also $J'_k=J_l\cap J_m$ for $\{k,l,m\}=\{1,2,3\}$.
It follows from \cite[Theorem 6.1]{cp:minsl} that, if $\lambda(h_{i_0})\ne 0$ and $\lambda$ is supported on the three connected components of $I\backslash \{i_0\}$, then $V_q(\lambda)$ has exactly three equivalence classes of minimal affinizations. Moreover:

\begin{thm}\label{t:regmma}
Let $\gb\lambda\in\cal P_q^+$ be such that $\wt(\gb\lambda)=\lambda$ where $\lambda$ is as above. Then $V_q(\gb\lambda)$ is a minimal affinization of $V_q(\gb\lambda)$ iff there exists $k\in\{1,2,3\}$ such that $V_q(\gb\lambda_{J_l})$ is a minimal affinization of $V_q(\lambda_{J_l})$ for $l\ne k$.\hfill\qedsymbol
\end{thm}

\begin{defn}
Given $\lambda\in P^+$ and $k\in\{1,2,3\}$, let $M_k(\lambda)$ be the quotient of $A(\lambda)$ by the submodule generated by the vectors $x_{\alpha,1}^-v_\lambda$ for all $\alpha\in R_{J_l}^+$ with $l\ne k$. Suppose $\gb\lambda\in\cal P_\mathbb A^{++}$ is such that $V_q(\gb\lambda_{J_l})$ is a minimal affinization of $V_q(\lambda_{J_l})$ for $l\ne k$. Set $T_k(\lambda)$ to be the $\lie g[t]$-submodule of $M(\lambda^{J'_k})\otimes L(\gb\lambda^{I\backslash J'_k})$ generated by the top weight space.
\end{defn}

It is quite simple to see that $M_k(\lambda)$ is a restricted $\lie g[t]$-module and that $M(\lambda)$ is a quotient of $M_k(\lambda)$ for all $k$. Moreover, proceeding similarly to the proofs of Propositions \ref{p:TqM} and \ref{p:MqN} we get the following analogue (we omit the details).

\begin{prop}
Let $\gb\lambda\in\cal P_\mathbb A^{++}$ and $k\in\{1,2,3\}$ be such that $V_q(\gb\lambda_{J_l})$ is a minimal affinization of $V_q(\lambda_{J_l})$ for $l\ne k$. Then there exist surjective $\lie g[t]$-module maps $M_k(\lambda)\twoheadrightarrow L(\gb\lambda)\twoheadrightarrow T_k(\lambda)$.\hfill\qedsymbol
\end{prop}

\begin{con}\label{cj:3a}
Suppose $\lambda\in P^+$ is supported on the three connected components of $I\backslash \{i_0\}$. Then, $T_k(\lambda)$ and $M_k(\lambda)$ are isomorphic for every $k\in\{1,2,3\}$.
\end{con}

\begin{cor}
Suppose $\gb\lambda\in\cal P_\mathbb A^{++}$ and $k\in\{1,2,3\}$ are such that $V_q(\gb\lambda_{J_l})$ is a minimal affinization of $V_q(\lambda_{J_l})$ for $l\ne k$ and $\wt(\gb\lambda)$ is supported on the three connected components of $I\backslash \{i_0\}$. Then, $T_k(\lambda)\cong L(\gb\lambda)\cong M_k(\lambda)$.
\end{cor}

We now prove Conjecture \ref{cj:3a} for $\lie g$ of type $D_4$. Thus, let $\lambda\in P^+$ be such that $\lambda(h_i)\ne 0$ for all $i\ne 2$, and let $\gb\lambda$ be such that $V_q(\gb\lambda)$ is a minimal affinization of $V_q(\lambda)$. Set also $J_1=\{1,2,3\}, J_2=\{1,2,4\}$, and $J_3=\{2,3,4\}$. Without loss of generality we can assume that $V_q(\gb\lambda_{J_1})$ and $V_q(\gb\lambda_{J_2})$ are minimal affinizations. We want to show that $T_3(\lambda)\cong M_3(\lambda)$ in this case. We  also assume that
\begin{equation}
\gb\lambda = \gb\omega_{1,a,\lambda(h_1)}\ \gb\omega_{2,aq^{\lambda(h_1)+\lambda(h_2)+1},\lambda(h_2)}\ \gb\omega_{3,aq^{\lambda(h_2)+\lambda(h_3)+1},\lambda(h_3)}\ \gb\omega_{4,aq^{\lambda(h_2)+\lambda(h_4)+1},\lambda(h_4)}
\end{equation}
for some $a\in\mathbb C^\times$. The case
\begin{equation*}
\gb\lambda = \gb\omega_{1,a,\lambda(h_1)}\ \gb\omega_{2,aq^{-(\lambda(h_1)+\lambda(h_2)+1)},\lambda(h_2)}\ \gb\omega_{3,aq^{-(\lambda(h_2)+\lambda(h_3)+1)},\lambda(h_3)}\ \gb\omega_{4,aq^{-(\lambda(h_2)+\lambda(h_4)+1)},\lambda(h_4)}
\end{equation*}
is proved similarly. If $\lambda(h_2)\ne 0$, these two cases cover all minimal affinizations such that $V_q(\gb\lambda_{J_1})$ and $V_q(\gb\lambda_{J_2})$ are also minimal affinizations. Otherwise, there are two more possibilities for $\gb\lambda$ (see the closing remark of subsection \ref{ss:multmai}).

Let $v$ be the image of $1$ in $M_3(\lambda)$. By the very definition of $M_3(\lambda)$ we have the following relations
\begin{gather}\label{e:n31}
x_{\alpha_i,1}^-v=x_{\alpha_2+\alpha_j,1}^-v=x_{\alpha_1+\alpha_2+\alpha_3,1}^-v=x_{\alpha_1+\alpha_2+\alpha_4,1}^-v=0
\end{gather}
for all $i,j\in I, j\ne 2$. Using the commutation relations $[x_\alpha^-,x_\beta^-]=x_{\alpha+\beta}^-$ (up to multiple) we also get
\begin{gather}\label{e:n32}
x_{\alpha,2}^-v=0 \qquad\forall\ \alpha\in R^+.
\end{gather}
Let $\vartheta_1=\sum_{i=1}^4\alpha_i = \omega_1+\omega_3+\omega_4-\omega_2, \vartheta_2=\vartheta_1-\alpha_1= \omega_3+\omega_4-\omega_1$, and $\theta=\theta_{1,1}=\vartheta_1+\alpha_2=\omega_2$. It follows that
\begin{equation*}
M_3(\lambda) = \sum_{\gbr r\in\mathbb Z_{\ge 0}^3}U(\lie n^-)\gbr y_\gbr rv
\end{equation*}
where
\begin{equation*}
\gbr y_\gbr r = (x_{\theta,1}^-)^{r_3}(x_{\vartheta_2,1}^-)^{r_2}(x_{\vartheta_1,1}^-)^{r_1}.
\end{equation*}
Since  $\{\vartheta_1,\vartheta_2,\theta\}$ is a linearly independent subset of $\lie h^*$, it follows, as before, that $m_\mu(M_3(\lambda))\le 1$ for every $\mu\in P^+$ and equality may occur only if $\mu = \lambda-r_1\vartheta_1-r_2\vartheta_2-r_3\theta$ for some $r_j\in\mathbb Z_{\ge 0}$. But such elements are dominant iff \begin{equation*}
r_1\le\lambda(h_1)+r_2, \qquad r_3\le\lambda(h_2)+r_1,\quad\text{and}\quad r_1+r_2\le \min\{\lambda(h_3),\lambda(h_4)\}.
\end{equation*}
Set
\begin{equation}\label{e:d4pset}
\cal D_3(\lambda) = \left\{\gbr r\in\mathbb Z_{\ge 0}^3: r_1\le\lambda(h_1),\ r_3\le\lambda(h_2),\ r_1+r_2\le \min\{\lambda(h_3),\lambda(h_4)\}\right\}.
\end{equation}
Proceeding similarly to the proof of Lemma \ref{l:ubb3} one proves that
\begin{equation}
M_3(\lambda) = \sum_{\gbr r\in\cal D_3(\lambda)}U(\lie n^-)\gbr y_\gbr rv.
\end{equation}
Given $\gbr r\in\mathbb Z_{\ge 0}^3$, define
\begin{equation}
\wt(\gbr r) = r_1\vartheta_1+r_2\vartheta_{2}+r_3\theta\qquad\text{and}\qquad \gr(\gbr r) = r_1+r_2+r_3.
\end{equation}
Since $\{\vartheta_1,\vartheta_2,\theta\}$ is linearly independent, it follows that $\wt$ is an injective function.
In order to complete the proof of Conjecture \ref{cj:3a} in this case, it suffices to prove that $m_\mu(T_3(\lambda))\ge 1$ if $\mu=\lambda-\wt(\gbr r)$ for some $\gbr r\in\cal D_3(\lambda)$. In particular, it will follow that
\begin{equation}\label{e:d4m}
M_3(\lambda)[l] = \opl_{\gbr r\in\cal D_3(\lambda):\gr(\gbr r)=l}^{} V(\lambda-\wt(\gbr r)).
\end{equation}
We begin by proving the following proposition.

\begin{prop}\label{p:nominspin}
Let $b\in \mathbb A^\times, \mu=m_3\omega_3+m_4\omega_4\in P^+$, and $\gb\mu = \gb\omega_{3,b,m_3}\gb\omega_{4,bq^{m_4-m_3},m_4}$. Then, $L(\gb\mu)[l] \cong V(\mu-l\vartheta_2)$ for $0\le l\le \min\{m_3,m_4\}$ and $L(\gb\mu)[l] =0$ otherwise.
\end{prop}

\begin{proof}
Let $v$ be a highest-weight vector of $L(\gb\mu)$. Quite clearly $v$ satisfies relations \eqref{e:n31} and \eqref{e:n32}. Moreover, proceeding as above, we get
\begin{equation*}
L(\gb\mu) = \sum_{\gbr r\in\cal D_3(\mu)}U(\lie g)\gbr y_\gbr rv = \opl_{r=0}^{\min\{m_3,m_4\}} U(\lie n^-)(x_{\vartheta_2,1}^-)^rv.
\end{equation*}
and, by Lemma \ref{l:heis} once more,
\begin{equation}\label{e:nominspin}
(x_{\vartheta_1,1}^-)^{r_1}(x_{\vartheta_2,1}^-)^{r_2}v = (x_{\alpha_1}^-)^{r_1}(x_{\vartheta_2,1}^-)^{r_1+r_2}v.
\end{equation}

Without loss of generality, assume $m_4\ge m_3\ge 1$ and observe that
$$\gb\mu = \left(\prod_{j=0}^{m_3-1} \gb\omega_{\omega_3+\omega_4,bq^{1-m_3+2j}}\right)\left(\prod_{j=0}^{m_4-m_3-1} \gb\omega_{4,bq^{m_3+1+2j}}\right).$$
Then, by Proposition \ref{p:cyclic} and its corollary, $V_q(\gb\mu)$ is the $U_q(\tlie g)$-submodule of $$\left(V_q(\gb\omega_{\omega_3+\omega_4,bq^{1-m_3}})\otimes\cdots\otimes V_q(\gb\omega_{\omega_3+\omega_4,bq^{m_3-1}}) \right) \left(V_q(\gb\omega_{4,bq^{m_3+1}})\otimes\cdots\otimes V_q(\gb\omega_{4,bq^{2m_4-m_3-1}}) \right)$$
generated by the top weight space. Let $M'(\omega_3+\omega_4)$ be the pullback of $\overline{V_q(\gb\omega_{\omega_3+\omega_4,bq^{m}})}$ by $\tau_b$, where $m\in\mathbb Z$ and let $T'(\mu)$ be the $\lie g[t]$-submodule of $M'(\omega_3+\omega_4)^{\otimes m_3}\otimes M(\omega_4)^{\otimes m_4-m_3-1}$. As before, it follows from Lemma \ref{l:morlim} that $T'(\mu)$ is a quotient of $L(\gb\mu)$. Hence, we are left to show that $T'(\mu)$ has $V(\mu-l\vartheta_2)$ as an irreducible $\lie g$-submodule for every $0\le l\le m_3$. Moreover, it suffices to consider the case $m_4=m_3=m\in\mathbb Z_{>0}$.
Observe that $V_q(\gb\omega_{\omega_3+\omega_4,bq^{m}})$ is not a minimal affinization and that $V(\omega_3)\otimes V(\omega_4)\cong V(\omega_3+\omega_4)\oplus V(\omega_1)$. In other words, the proposition is proved for $m_3=m_4=1$. Finally, let $v_j, j=1,\dots,m$, be a highest-weight vector of the $j$-th copy of $M'(\omega_3+\omega_4)$ in $M'(\omega_3+\omega_4)^{\otimes m}$ and let $v_j^0$ be its image in $M'(\omega_3+\omega_4)(0)$. Then,
$$(x_{\vartheta_2,1}^-)^l(v_1\otimes v_2\otimes\cdots\otimes v_l\otimes v_{l+1}^0\otimes \cdots\otimes v_m^0) = (x_{\vartheta_2,1}^-v_1)\otimes (x_{\vartheta_2,1}^-v_2)\otimes\cdots\otimes (x_{\vartheta_2,1}^-v_l)\otimes v_{l+1}^0\otimes \cdots\otimes v_m^0$$
and we are done using a simple induction on $l$.
\end{proof}

Let $v_1$ be a highest-weight vector of $M(\lambda^{\{1,2\}})$ and $v_2$ be a highest-weight vector of $L(\gb\lambda^{\{3,4\}})$. It follows from Proposition \ref{p:nd3} and \eqref{e:nominspin} that, if $\gbr r\in\cal D_3(\lambda)$, then
\begin{equation}
\gbr y_\gbr r (v_1\otimes v_2) = ((x_{\theta,1}^-)^{r_3}v_1)\otimes ((x_{\alpha_1}^-)^{r_1}(x_{\vartheta_2,1}^-)^{r_1+r_2}v_2).
\end{equation}
The proof of \eqref{e:d4m} is now completed similarly to the end of the proof of Proposition \ref{p:Nb3}.

\subsection{Multiple minimal affinizations: the irregular case}\label{ss:multmai}
Keep the notation of subsection \ref{ss:multmar}. If $\lambda$ is supported on the three connected components of $I\backslash \{i_0\}$ and $\lambda(h_{i_0})= 0$, it follows from \cite{cp:minirr} that the number of equivalence classes of minimal affinizations of $V_q(\lambda)$ is not uniformly bounded (it grows as $\lambda$ ``grows''). If $\lie g$ is of type $E$, write $I$ as the disjoint union of two connected subdiagrams of type $A$, say $I_1$ and $I_2$, and the subdiagram of type $D_4$, say $J$. For $\lie g$  of type $D$ we write $I$ as the disjoint union of a subdiagram $I_1$ of type $A$ and the subdiagram $J$ of type $D_4$ (for convenience we set $I_2=\emptyset$ and $\lambda^{\emptyset}=0$). Similarly to the proof of Proposition \ref{p:TqM} one proves:

\begin{prop}
Let $\lambda\in P^+$ and $\gb\lambda\in\cal P_q^+$ be such that $V_q(\gb\lambda)$ is a minimal affinization of $V_q(\lambda)$. Then $L(\gb\lambda)$ projects onto the $\lie g[t]$-submodule of $L(\lambda^{I_1})\otimes L(\gb\lambda^J)\otimes L(\lambda^{I_2})$ generated by the top weight space.\hfill\qedsymbol
\end{prop}

The natural conjecture is then stated as:

\begin{con}
Let $\lambda\in P^+$ and $\gb\lambda\in\cal P_q^+$ be such that $V_q(\gb\lambda)$ is a minimal affinization of $V_q(\lambda)$. Then $L(\gb\lambda)$ is isomorphic to the $\lie g[t]$-submodule of $L(\lambda^{I_1})\otimes L(\gb\lambda^J)\otimes L(\lambda^{I_2})$ generated by the top weight space.
\end{con}

We shall leave this conjecture in a purely speculative tone for the moment and postpone further discussion of these cases to a forthcoming publication.

\begin{rem}
Let $\gb\lambda$ be as in \cite[Theorem 2.2 (a)$_{3,4}$ or (b)$_{3,4}$]{cp:minirr}.
If conditions (a)$_{3,4}$ are satisfied, the results of subsection \ref{ss:multmar} apply and, hence, the graded character of $L(\gb\lambda)$ is given by the right-hand-side of \eqref{e:d4m}. In order to prove the conjecture of remark (1) that follows Theorem 2.2 of \cite{cp:minirr}, it suffices to show that, if  conditions (b)$_{3,4}$ are satisfied, then the graded character of $L(\gb\lambda)$ is also given by the right-hand-side of \eqref{e:d4m}. The proof is essentially the same as for the former case replacing Proposition \ref{p:nominspin} by its appropriate obvious modification. We omit the details.
\end{rem}

\bibliographystyle{amsplain}

\end{document}